# Deadzone-Adapted Disturbance Suppression Control for Strict-Feedback Systems

**Iasson Karafyllis[*], Miroslav Krstic[**] and Alexandros Aslanidis[*]**

[*]Dept. of Mathematics, National Technical University of Athens, Zografou Campus, 15780, Athens, Greece, email: iasonkar@central.ntua.gr; iasonkaraf@gmail.com

[**]Dept. of Mechanical and Aerospace Eng., University of California, San Diego, La Jolla, CA 92093-0411, U.S.A., email: krstic@ucsd.edu

**Abstract**

In this paper we extend our recently proposed Deadzone-Adapted Disturbance Suppression (DADS) Control approach from systems with matched uncertainties to general systems in parametric strict feedback form. The DADS approach prevents gain and state drift regardless of the size of the disturbance and unknown parameter and achieves an attenuation of the plant output to an assignable small level, despite the presence of persistent disturbances and unknown parameters of arbitrary and unknown bounds. The controller is designed by means of a step-by-step backstepping procedure which can be applied in an algorithmic fashion. Examples are provided which illustrate the efficiency of the DADS controller compared to existing adaptive control schemes.

**Keywords:** Robust Adaptive Control, Feedback Stabilization, Leakage.

## 1. Introduction

The main objective of adaptive control is the regulation of the plant state. Various control strategies have been proposed in the literature, addressing this problem: indirect adaptive control schemes rely on an adjustable estimate of the unknown parameters (see for instance [10, 14, 15, 16, 17, 18]) whereas direct control schemes adjust the parameters of the controller itself (see for instance [1,12]). To guarantee robustness with respect to external disturbances for time-invariant nonlinear control systems for which no persistence of excitation condition is assumed, several strategies have been proposed: nonlinear damping (see [9, 10, 15, 17]), leakage (see [2, 3, 22, 25]), projection methodologies (see [2] and Appendix E in [15]), supervision for direct adaptive schemes (see [1]), dynamic (high) gains or gain adjustment (see [4, 5, 6, 13, 21]) and deadzone in the update law (introduced in [23] and well explained in the book [2]). Some of these control schemes require no knowledge of the size of the disturbance (nonlinear damping), while others (e.g., projection) do require knowledge of the size of the unknown parameters. Schemes like $\sigma$-modification (leakage) regulate the plant state to a neighborhood of zero whose size depends on the unknown parameter.

Recently, we introduced in [12] a novel direct adaptive control approach, the Deadzone-Adapted Disturbance Suppression (DADS), for a class of systems that satisfy the matching condition (see [15]). The controller proposed in [12], under no assumptions on the unknown parameters, achieves regulation of the plant state to a ball around the origin of assignable radius, even in the presence of persistent ($L^\infty$) disturbances. The DADS controller encompasses dynamic nonlinear damping terms involving a single adjustable gain and a deadzone in the update law of the gain. Moreover, it was shown in [12], that there are systems, not satisfying the matching condition, which do not admit a feedback law (dynamic or static) that guarantees attenuation of the plant state to an assignable level despite the presence of disturbances and unknown parameters of arbitrary and unknown bounds.

The present paper is devoted to the advancement of the DADS approach to systems in parametric strict feedback form. Inspired by the idea of dynamic gain used in [4, 5, 6, 13, 19, 20, 21, 24] and the recent work [10] where indirect, adaptive feedback designs which guarantee KL estimates in the disturbance-free case were provided, here we extend the DADS approach to the triangular nonlinear case where no matching condition holds. The DADS controller prevents gain and state drift regardless of the size of the disturbance and unknown parameter and achieves an attenuation of a plant output to an assignable small level, despite the presence of persistent disturbances and unknown parameters of arbitrary and unknown bounds. The controller is designed by means of a step-by-step backstepping procedure (see Lemma 4 below) which can be applied in an algorithmic fashion similar to the backstepping procedures that are provided in [10, 14, 15, 16, 17, 18].

The combination of properties that the DADS controller achieves is unique in the sense that no other adaptive control scheme in the literature can guarantee these specific properties for systems for which no persistence of excitation condition is assumed and no bound for the disturbance and the unknown parameters is known. An exception is the adaptive control scheme proposed in [11] with delays, which however has been applied only to a very simple (scalar) system. The DADS controller achieves almost all properties that the delay scheme in [11] can achieve without using delays. Therefore, the DADS controller is easier to implement.

The paper is organized as follows. In Section 2, we present the controller and give the main stability results. In Section 3, we give an illustrative example comparing the DADS controller with an adaptive controller based on $\sigma$-modification. Section 4 provides the proofs of all results while the concluding remarks of the present work are given in Section 5.

**Notation and Basic Notions.** Throughout this paper, we adopt the following notation.
* $\mathbb{R}_+ := [0, +\infty)$. For a vector $x \in \mathbb{R}^n$, $|x|$ denotes its Euclidean norm and $x'$ denotes its transpose. We use the notation $x^+$ for the positive part of the real number $x \in \mathbb{R}$, i.e., $x^+ = \max(x, 0)$. For a matrix $A \in \mathbb{R}^{n \times m}$, $A' \in \mathbb{R}^{m \times n}$ denotes its transpose. We write $A \leq B$ for two symmetric matrices $A, B \in \mathbb{R}^{n \times n}$ when the matrix $(B - A) \in \mathbb{R}^{n \times n}$ is positive definite, i.e., when $x'(B-A)x > 0$ for all $x \in \mathbb{R}^n \setminus \{0\}$.

* Let $D \subseteq \mathbb{R}^n$ be an open set and let $S \subseteq \mathbb{R}^n$ be a set that satisfies $D \subseteq S \subseteq cl(D)$, where $cl(D)$ is the closure of $D$. By $C^0(S; \Omega)$, we denote the class of continuous functions on $S$, which take values in $\Omega \subseteq \mathbb{R}^m$. By $C^k(S; \Omega)$, where $k \geq 1$ is an integer, we denote the class of functions on $S \subseteq \mathbb{R}^n$, which take values in $\Omega \subseteq \mathbb{R}^m$ and have continuous derivatives of order $k$. In other words, the functions of class $C^k(S; \Omega)$ are the functions which have continuous derivatives of order $k$ in $D = \text{int}(S)$ that can be continued continuously to all points in $\partial D \cap S$. When $\Omega = \mathbb{R}$ then we write $C^0(S)$ or $C^k(S)$. A function $f \in C^\infty(S; \Omega) = \bigcap_{k=0}^{\infty} C^k(S; \Omega)$ is called a smooth function.

* Let the non-empty set $\Theta \subseteq \mathbb{R}^p$ be given. By $L^\infty(\mathbb{R}_+;\Theta)$ we denote the class of essentially bounded, Lebesgue measurable functions $d:\mathbb{R}_+ \to \Theta$. For $d \in L^\infty(\mathbb{R}_+;\Theta)$ we define $\|d\|_\infty = \sup_{t \geq 0}(|d(t)|)$, where $\sup_{t \geq 0}(|d(t)|)$ is the essential supremum.

* By $K$ we denote the class of increasing continuous functions $a:\mathbb{R}_+ \to \mathbb{R}_+$ with $a(0) = 0$. By $K_\infty$ we denote the class of increasing continuous functions $a:\mathbb{R}_+ \to \mathbb{R}_+$ with $a(0) = 0$ and $\lim_{s \to +\infty}(a(s)) = +\infty$. By $KL$ we denote the set of all continuous functions $\sigma:\mathbb{R}_+ \times \mathbb{R}_+ \to \mathbb{R}_+$ with the properties: (i) for each $t \geq 0$ the mapping $\sigma(\cdot,t)$ is of class $K$; (ii) for each $s \geq 0$, the mapping $\sigma(s,\cdot)$ is non-increasing with $\lim_{t \to +\infty}(\sigma(s,t)) = 0$.

We next recall certain standard notions in the analysis of output stability: see for instance [8, 9, 26, 27].

Let $f:\mathbb{R}^n \to \mathbb{R}^n$ be a locally Lipschitz vector field with $f(0) = 0$ and $h:\mathbb{R}^n \to \mathbb{R}^p$ be a continuous mapping with $h(0) = 0$. Consider the dynamical system

$$\dot{x} = f(x), \; x \in \mathbb{R}^n \tag{1.1}$$

with output

$$Y = h(x) \tag{1.2}$$

We assume that the dynamical system (1.1) is forward complete, i.e., for every $x_0 \in \mathbb{R}^n$ the unique solution $x(t) = \phi(t,x_0)$ of the initial-value problem (1.1) with initial condition $x(0) = x_0$ exists for all $t \geq 0$. We use the notation $Y(t,x_0) = h(\phi(t,x_0))$ for all $t \geq 0$, $x_0 \in \mathbb{R}^n$ and $B_R = \{x \in \mathbb{R}^n : |x| < R\}$ for all $R > 0$. We say that system (1.1), (1.2) is

i) *Lagrange output stable* if for every $R > 0$ the set $\{|Y(t,x_0)| : x_0 \in B_R, t \geq 0\}$ is bounded.

ii) *Lyapunov output stable* if for every $\varepsilon > 0$ there exists $\delta(\varepsilon) > 0$ such that for all $x_0 \in B_{\delta(\varepsilon)}$, it holds that $|Y(t,x_0)| \leq \varepsilon$ for all $t \geq 0$.

iii) *Globally Asymptotically Output Stable (GAOS)* if system (1.1), (1.2) is Lagrange and Lyapunov output stable and $\lim_{t \to +\infty}(Y(t,x_0)) = 0$ for all $x_0 \in \mathbb{R}^n$.

iv) *Uniformly Globally Asymptotically Output Stable (UGAOS)* if system (1.1), (1.2) is Lagrange and Lyapunov output stable and for every $\varepsilon, R > 0$ there exists $T(\varepsilon, R) > 0$ such that for all $x_0 \in B_R$, it holds that $|Y(t,x_0)| \leq \varepsilon$ for all $t \geq T(\varepsilon, R)$.

It should be noted that (see Theorem 2.2 on page 62 in [8]) UGAOS for system (1.1), (1.2) is equivalent to the existence of a function $\beta \in KL$ such that the following estimate holds for all $x_0 \in \mathbb{R}^n$ and $t \geq 0$:

$$|Y(t,x_0)| \leq \beta(|x_0|,t) \tag{1.3}$$

We say that system (1.1), (1.2) is *practically Uniformly Globally Asymptotically Output Stable (p-UGAOS)* if there exists a function $\beta \in KL$ and a constant $\alpha > 0$ such that the following estimate holds for all $x_0 \in \mathbb{R}^n$ and $t \geq 0$:

$$|Y(t,x_0)| \leq \beta(|x_0|,t) + \alpha \tag{1.4}$$

When $h(x) = x$ then the word "output" in the above properties is omitted (e.g., Lagrange stability, Lyapunov stability, GAS, UGAS, p-UGAS).

Let $f : \mathbb{R}^n \times \mathbb{R}^p \to \mathbb{R}^n$ be a locally Lipschitz mapping with $f(0,0) = 0$. Consider the control system

$$\dot{x} = f(x,d), \ x \in \mathbb{R}^n, \ d \in \mathbb{R}^p \tag{1.5}$$

We assume that system (1.5) is forward complete, i.e., for every $x_0 \in \mathbb{R}^n$ and for every Lebesgue measurable and locally essentially bounded input $d : \mathbb{R}_+ \to \mathbb{R}^p$ the unique solution $x(t) = \phi(t, x_0; d)$ of the initial-value problem (1.5) with initial condition $x(0) = x_0$ corresponding to input $d : \mathbb{R}_+ \to \mathbb{R}^p$ exists for all $t \geq 0$. We use the notation $Y(t, x_0; d) = h(\phi(t, x_0; d))$ for all $t \geq 0$, $x_0 \in \mathbb{R}^n$ and for every Lebesgue measurable and locally essentially bounded input $d : \mathbb{R}_+ \to \mathbb{R}^p$.

We say that system (1.5), (1.2) is *Input-to-Output Stable (IOS)* if there exist functions $\beta \in KL$, $\gamma \in K$ such that the following estimate holds for all $x_0 \in \mathbb{R}^n$, $t \geq 0$ and for every $d \in L^\infty(\mathbb{R}_+)$:

$$|Y(t, x_0; d)| \leq \beta(|x_0|, t) + \gamma(\|d\|_\infty) \tag{1.6}$$

We say that system (1.5), (1.2) is *practically Input-to-Output Stable (p-IOS)* if there exist functions $\beta \in KL$, $\gamma \in K$ and a constant $\alpha > 0$ such that the following estimate holds for all $x_0 \in \mathbb{R}^n$, $t \geq 0$ and for every $d \in L^\infty(\mathbb{R}_+)$:

$$|Y(t, x_0; d)| \leq \beta(|x_0|, t) + \gamma(\|d\|_\infty) + \alpha \tag{1.7}$$

We say that system (1.5), (1.2) satisfies the *practical Output Asymptotic Gain (p-OAG)* property if there exists a non-decreasing continuous function $\gamma : \mathbb{R}_+ \to \mathbb{R}_+$ with $\gamma(0) = 0$ and a constant $\alpha > 0$ such that the following estimate holds for all $x_0 \in \mathbb{R}^n$ and for every $d \in L^\infty(\mathbb{R}_+)$:

$$\limsup_{t \to +\infty} \left(|Y(t, x_0; d)|\right) \leq \gamma(\|d\|_\infty) + \alpha \tag{1.8}$$

When $\gamma \equiv 0$ we say that system (1.5), (1.2) satisfies the *zero practical Output Asymptotic Gain property (zero p-OAG)*.

When $h(x) = x$ then the word "output" in the above properties is either replaced by the word "state" (e.g., ISS, p-ISS) or is omitted (e.g., p-AG, zero p-AG).

We say that system (1.5), (1.2) satisfies the *practical Uniform Bounded-Input-Bounded-State (p-UBIBS)* property if there exists a function $\gamma \in K_\infty$ and a constant $\alpha > 0$ such that the following estimate holds for all $x_0 \in \mathbb{R}^n$ and for every $d \in L^\infty(\mathbb{R}_+)$:

$$\sup_{t \geq 0}\left(|\phi(t, x_0; d)|\right) \leq \gamma(|x_0|) + \gamma(\|d\|_\infty) + \alpha \tag{1.9}$$

Clearly, the p-UBIBS property is equivalent to the existence of a continuous function $B : \mathbb{R}^n \times \mathbb{R}_+ \to \mathbb{R}_+$ for which the following estimate holds for all $x_0 \in \mathbb{R}^n$ and for every $d \in L^\infty(\mathbb{R}_+)$:

$$\sup_{t \geq 0}\left(|\phi(t, x_0; d)|\right) \leq B(x_0, \|d\|_\infty) \tag{1.10}$$

In the disturbance-free case, the p-UBIBS property coincides with the Lagrange stability property.

## 2. Main Results

In this work we study nonlinear control systems of the form

$$\dot{x}_i = x_{i+1}, \text{ for } i = 1,...,n \tag{2.1}$$

$$\dot{y}_j = h_j(x, y_1,..., y_j) + g_j(x, y_1,..., y_j, \theta) y_{j+1} + \varphi'_j(x, y_1,..., y_j)\theta + \alpha'_j(x, y_1,..., y_j)d,$$
$$\text{for } j = 1,...,m \tag{2.2}$$

where $x_{n+1} = y_1$, $y_{m+1} = u$, $(x, y) \in \mathbb{R}^n \times \mathbb{R}^m$ is the state with $x = (x_1,..., x_n)$, $y = (y_1,..., y_m)$, $u \in \mathbb{R}$ is the control input, $d \in \mathbb{R}^l$, $\theta \in \Theta \subseteq \mathbb{R}^p$ are disturbances, $h_j : \mathbb{R}^n \times \mathbb{R}^j \to \mathbb{R}$, $h_j : \mathbb{R}^n \times \mathbb{R}^j \times \Theta \to \mathbb{R}$, $\varphi_j : \mathbb{R}^n \times \mathbb{R}^j \to \mathbb{R}^p$, $\alpha_j : \mathbb{R}^n \times \mathbb{R}^j \to \mathbb{R}^l$ for $j = 1,...,m$ are smooth mappings with $h_j(0,0) = 0$ and $\varphi_j(0,0) = 0$ for $j = 1,...,m$. Systems of the form (2.1), (2.2) are termed in the literature as triangular systems or systems in parametric strict feedback form. It should be noted that $\theta \in \Theta \subseteq \mathbb{R}^p$ is a vanishing perturbation while $d \in \mathbb{R}^l$ is -in general- a non-vanishing perturbation.

We assume that there exist smooth positive mappings $\eta_j : \mathbb{R}^n \times \mathbb{R}^j \to (0, +\infty)$ for $j = 1,...,m$, $\mu_j : \mathbb{R}^n \times \mathbb{R}^j \to (0, +\infty)$ for $j = 1,...,m-1$ such that the following inequalities hold:

$$\eta_j(x, y_1,..., y_j) \le g_j(x, y_1,..., y_j, \theta),$$
$$\text{for all } (x, y_1,..., y_j) \in \mathbb{R}^n \times \mathbb{R}^j, \theta \in \Theta, j = 1,...,m \tag{2.3}$$

$$g_j(x, y_1,..., y_j, \theta) \le \mu_j(x, y_1,..., y_j)(1+|\theta|),$$
$$\text{for all } (x, y_1,..., y_j) \in \mathbb{R}^n \times \mathbb{R}^j, \theta \in \Theta, j = 1,...,m-1 \tag{2.4}$$

When $\Theta \subseteq \mathbb{R}^p$ is a bounded set, we can exploit (2.3) and (2.4) in order to design a smooth feedback law $u = k(x, y)$ that guarantees the ISS property with respect to $d \in \mathbb{R}^l$ (uniformly in $\theta \in \Theta$); see [8]. On the other hand, when $\theta \in \Theta$ is constant (but unknown) and $d \equiv 0$ then we can apply the methodologies in [15] in order to design an adaptive controller that guarantees GAOS for the closed-loop system with output $Y = (x, y)$ (notice that here we are talking about GAOS and not GAS since the closed-loop system has additional states: the states of the controller). When $\theta \in \Theta$ is constant (and unknown) and $d$ is not necessarily zero one can employ standard techniques in order to achieve the p-ISS property with respect to $d \in \mathbb{R}^l$ for the closed-loop system. One such technique that does not require the knowledge of bounds for $\theta \in \Theta$ is the so-called $\sigma$ – modification which guarantees the p-ISS property for the closed-loop system. However, for the $\sigma$ – modification the radius of the residual set is an increasing function of the unknown $|\theta|$, which means that the output $Y = (x, y)$ (even in the case $d \equiv 0$) approaches a neighborhood of zero of unknown (and possibly large) radius.

An alternative approach to $\sigma$ – modification was proposed in [12] for systems that satisfy a matching condition: the so-called Deadzone-Adapted Disturbance Suppression (DADS) controller. The DADS controller combines both the adaptation idea and the robust control methodology and uses only one integrator with deadzone. The DADS controller guarantees an attenuation of the plant output to an assignable small level, despite the presence of disturbances and unknown parameters of arbitrary and unknown bounds while in the same time prevents gain and state drift regardless of the size of the disturbance and unknown parameter. In this work, we extend the DADS controller to the case (2.1), (2.2) where no matching condition holds. Our main result is stated next.

**Theorem 1 (DADS controller for systems in parametric strict feedback form):** *Consider the control system (2.1), (2.2) and assume that there exist smooth positive mappings $\eta_j : \mathbb{R}^n \times \mathbb{R}^j \to (0, +\infty)$ for $j = 1, ..., m$, $\mu_j : \mathbb{R}^n \times \mathbb{R}^j \to (0, +\infty)$ for $j = 1, ..., m-1$ such that (2.3) and (2.4) hold. Then for every constants $b, \Gamma, \varepsilon, c, a > 0$ and for every $\kappa, \lambda \in K_\infty \cap C^\infty(\mathbb{R}_+)$ there exist functions $k, V \in C^\infty(\mathbb{R}^n \times \mathbb{R}^m \times \mathbb{R})$ with $V(0, 0, z) = k(0, 0, z) = 0$ for $z \in \mathbb{R}$ and $V(x, y, z) > 0$ for $(x, y) \neq 0$, a constant $M \geq 1$ and a function $B \in C^0(\mathbb{R}^n \times \mathbb{R}^m \times \mathbb{R} \times \mathbb{R}_+ \times \mathbb{R}_+)$ such that for every $(x_0, y_0, z_0) \in \mathbb{R}^n \times \mathbb{R}^m \times \mathbb{R}$, $d \in L^\infty(\mathbb{R}_+; \mathbb{R}^l)$, $\theta \in L^\infty(\mathbb{R}_+; \Theta)$ the unique solution of (2.1), (2.2) with*

$$u = k(x, y, z) \tag{2.5}$$

$$\dot{z} = \Gamma e^{-z} \left( V(x, y, z) - \frac{\varepsilon^2}{2M} \right)^+ \tag{2.6}$$

*and $(x(0), y(0), z(0)) = (x_0, y_0, z_0)$ that corresponds to $d \in L^\infty(\mathbb{R}_+; \mathbb{R}^l)$, $\theta \in L^\infty(\mathbb{R}_+; \Theta)$ exists for all $t \geq 0$ and satisfies the following estimates for all $t \geq 0$:*

$$\limsup_{t \to +\infty} (|x(t)| + |y_1(t)|) \leq \varepsilon, \quad \limsup_{t \to +\infty} (V(x(t), y(t), z(t))) \leq \frac{\varepsilon^2}{2M} \tag{2.7}$$

$$z_0 \leq z(t) \leq \lim_{s \to +\infty} (z(s)) \leq B(x_0, y_0, z_0, \|d\|_\infty, \|\theta\|_\infty) \tag{2.8}$$

$$|x(t)| + |y(t)| \leq B(x_0, y_0, z_0, \|d\|_\infty, \|\theta\|_\infty) \tag{2.9}$$

$$|x(t)|^2 + |y_1(t)|^2 \leq MV(x(t), y(t), z(t))$$
$$\leq Me^{-ct} V(x_0, y_0, z_0) + a \frac{\|d\|_\infty^2 + \left( \left( \|\theta\|_\infty - b - \lambda(e^{z_0}) \right)^+ \right)^2}{c(1 + \kappa(e^{z_0}))} \tag{2.10}$$

*Moreover, if $\lim_{t \to +\infty} (d(t)) = 0$ and $\limsup_{t \to +\infty} (|\theta(t)|) \leq b + \lambda\left( e^{\lim_{t \to +\infty} (z(t))} \right)$ then $\lim_{t \to +\infty} (|x(t)| + |y(t)|) = 0$.*

Inequality (2.7) guarantees that the plant output $Y = (x, y_1)$ is eventually led (even in the case where a disturbance is present) to a ball around the origin of <u>*assignable radius*</u>. Indeed, inequality (2.7) guarantees the zero p-OAG property $\limsup_{t \to +\infty} (|(x(t), y_1(t))|) \leq \alpha$ with a constant $\alpha > 0$ independent of $\theta$. This is a feature that cannot be guaranteed by $\sigma$-modification when $\theta$ is constant. Moreover, the state, i.e., both $(x(t), y(t))$ and $z(t)$, remain bounded for every initial condition and every disturbance $d \in L^\infty(\mathbb{R}_+)$-see (2.8), (2.9); this is the p-UBIBS property. Inequality (2.10) shows that the p-IOS property from the disturbance $d \in L^\infty(\mathbb{R}_+; \mathbb{R}^l)$ holds for the closed-loop system (2.1), (2.2), (2.5), (2.6) with output $Y = (x, y_1)$. When $d \equiv 0$ and $\theta$ is constant then the closed-loop system (2.1), (2.2), (2.5), (2.6) with output $Y = (x, y_1)$: i) is Lagrange stable, ii) satisfies the p-UGAOS property, and iii) has a globally attracting set, namely the set $\left\{ (x, y, z) \in \mathbb{R}^n \times \mathbb{R}^m \times \mathbb{R} : V(x, y, z) \leq \dfrac{\varepsilon^2}{2M} \right\}$. Finally, it should also be noticed that convergence to $(x, y) = (0, 0)$ is not excluded: it can be achieved when $\lim_{t \to +\infty}(d(t)) = 0$ and $\limsup_{t \to +\infty}(|\theta(t)|) \leq b + \lambda\left(e^{\lim_{t \to +\infty}(z(t))}\right)$. This is another feature that cannot be guaranteed by $\sigma-$modification when $\theta$ is constant. For the DADS controller in the disturbance-free case we have the following conclusion: *if the plant state is not led asymptotically to zero then the estimate* $\lambda\left(e^{z(t)}\right) \leq \lambda\left(e^{\lim_{t \to +\infty}(z(t))}\right) < |\theta| - b$ *holds for all* $t \geq 0$.

The proof of Theorem 1 is constructive and is based on a backstepping lemma (see Lemma 4 in Section 4). Moreover, Theorem 1 is based on the following stability result which is interesting on its own.

**Theorem 2 (Lyapunov conditions for systems with deadzone):** *Consider the system*

$$\dot{x} = f(x, z, \theta, d)$$
$$\dot{z} = \Gamma e^{-z}(V(x, z) - \varepsilon)^+ \tag{2.11}$$

*where* $(x, z) \in \mathbb{R}^n \times \mathbb{R}$ *is the state,* $\theta \in \Theta \subseteq \mathbb{R}^p$, $d \in D \subseteq \mathbb{R}^l$ *are disturbances,* $\Gamma, \varepsilon > 0$ *are constants,* $f : \mathbb{R}^n \times \mathbb{R} \times \Theta \times D \to \mathbb{R}^n$ *is a locally Lipschitz mapping with* $f(0, z, \theta, 0) = 0$ *for all* $(z, \theta) \in \mathbb{R} \times \Theta$ *and* $V \in C^1(\mathbb{R}^n \times \mathbb{R})$ *with* $V(0, z) = 0$ *for all* $z \in \mathbb{R}$ *and* $V(x, z) > 0$ *for all* $x \neq 0$, $z \in \mathbb{R}$ *is a function that satisfies the following properties:*

**(i)** *For every* $M > 0$ *the set* $\{(x, z) \in \mathbb{R}^n \times [-M, M] : V(x, z) \leq M\}$ *is bounded.*

**(ii)** *There exist functions* $\kappa, \lambda \in K_\infty$ *and constants* $a, b, c > 0$ *such that the following inequality holds for all* $(x, z) \in \mathbb{R}^n \times \mathbb{R}, \theta \in \Theta, d \in D$:

$$\dfrac{\partial V}{\partial x}(x, z) f(x, z, \theta, d) + \Gamma e^{-z} \dfrac{\partial V}{\partial z}(x, z)(V(x, z) - \varepsilon)^+$$
$$\leq -cV(x, z) + a \dfrac{|d|^2 + \left(\left(|\theta| - b - \lambda(e^z)\right)^+\right)^2}{1 + \kappa(e^z)} \tag{2.12}$$

*Then there exists a function* $B \in C^0(\mathbb{R}^n \times \mathbb{R} \times \mathbb{R}_+ \times \mathbb{R}_+)$ *such that for every* $(x_0, z_0) \in \mathbb{R}^n \times \mathbb{R}$, $d \in L^\infty(\mathbb{R}_+; D)$, $\theta \in L^\infty(\mathbb{R}_+; \Theta)$ *the unique solution of (2.11) with* $(x(0), z(0)) = (x_0, z_0)$ *that corresponds to* $d \in L^\infty(\mathbb{R}_+; D)$, $\theta \in L^\infty(\mathbb{R}_+; \Theta)$ *exists for all* $t \geq 0$ *and satisfies the following estimates:*

$$\limsup_{t \to +\infty}(V(x(t), z(t))) \leq \varepsilon \tag{2.13}$$

$$z_0 \leq z(t) \leq \lim_{s \to +\infty}(z(s)) \leq B\left(x_0, z_0, \|d\|_\infty, \|\theta\|_\infty\right), \text{ for all } t \geq 0 \tag{2.14}$$

$$|x(t)| \leq B\left(x_0, z_0, \|d\|_\infty, \|\theta\|_\infty\right), \text{ for all } t \geq 0 \tag{2.15}$$

$$V(x(t), z(t)) \leq e^{-ct}V(x_0, z_0) + a\frac{\|d\|_\infty^2 + \left(\left(\|\theta\|_\infty - b - \lambda\left(e^{z_0}\right)\right)^+\right)^2}{c\left(1 + \kappa\left(e^{z_0}\right)\right)}, \text{ for all } t \geq 0 \tag{2.16}$$

$$\limsup_{t \to +\infty}(V(x(t), z(t)))$$

$$\leq a\frac{\left(\limsup_{t \to +\infty}(|d(t)|)\right)^2 + \left(\left(\limsup_{t \to +\infty}(|\theta(t)|) - b - \lambda\left(e^{\lim_{t \to +\infty}(z(t))}\right)\right)^+\right)^2}{c\left(1 + \kappa\left(e^{\lim_{t \to +\infty}(z(t))}\right)\right)} \tag{2.17}$$

**Remark:** Estimates (2.14), (2.15) and (2.17) show that $\lim_{t \to +\infty}(|x(t)|) = 0$ provided that $\lim_{t \to +\infty}(d(t)) = 0$ and $\limsup_{t \to +\infty}(|\theta(t)|) \leq b + \lambda\left(e^{\lim_{t \to +\infty}(z(t))}\right)$.

Using the same tools we can also study nonlinear systems of the form

$$\dot{x}_i = h_i(x_1, ..., x_i) + g_i(x_1, ..., x_i, \theta)x_{i+1} + \varphi_i'(x_1, ..., x_i)\theta + \alpha_i'(x_1, ..., x_i)d, \text{ for } i = 1, ..., n \tag{2.18}$$

where $x_{n+1} = u$, $x = (x_1, ..., x_n) \in \mathbb{R}^n$ is the state, $u \in \mathbb{R}$ is the control input, $d \in \mathbb{R}^l$, $\theta \in \Theta \subseteq \mathbb{R}^p$ are disturbances, $h_i : \mathbb{R}^i \to \mathbb{R}$, $g_i : \mathbb{R}^i \times \Theta \to \mathbb{R}$, $\varphi_i : \mathbb{R}^i \to \mathbb{R}^p$, $\alpha_i : \mathbb{R}^i \to \mathbb{R}^l$ for $i = 1, ..., n$ are smooth mappings with $h_i(0,0) = 0$ and $\varphi_i(0,0) = 0$ for $i = 1, ..., n$. In this case we obtain the following result.

**Theorem 3 (DADS controller for systems in parametric strict feedback form with no integrators):** *Consider the control system (2.18) and assume that there exist smooth positive mappings $\eta_i : \mathbb{R}^i \to (0, +\infty)$ for $i = 1, ..., n$, $\mu_i : \mathbb{R}^i \to (0, +\infty)$ for $i = 1, ..., n-1$ such that the following inequalities hold:*

$$\eta_i(x_1, ..., x_i) \leq g_i(x_1, ..., x_i, \theta),$$
$$\text{for all } (x_1, ..., x_i) \in \mathbb{R}^i, \ \theta \in \Theta, \ i = 1, ..., n \tag{2.19}$$

$$g_i(x_1, ..., x_i, \theta) \leq \mu_i(x_1, ..., x_i)(1 + |\theta|),$$
$$\text{for all } (x_1, ..., x_i) \in \mathbb{R}^i, \ \theta \in \Theta, \ i = 1, ..., n-1 \tag{2.20}$$

*Then for every constants $b, \Gamma, \varepsilon, a, c > 0$ and for every $\kappa, \lambda \in K_\infty \cap C^\infty(\mathbb{R}_+)$ there exist functions $k, V \in C^\infty(\mathbb{R}^n \times \mathbb{R})$ with $V(0, z) = k(0, z) = 0$ for $z \in \mathbb{R}$ and $V(x, z) > 0$ for $x \neq 0$ and a function $B \in C^0(\mathbb{R}^n \times \mathbb{R} \times \mathbb{R}_+ \times \mathbb{R}_+)$ such that for every $(x_0, z_0) \in \mathbb{R}^n \times \mathbb{R}$, $d \in L^\infty(\mathbb{R}_+; \mathbb{R}^l)$, $\theta \in L^\infty(\mathbb{R}_+; \Theta)$ the unique solution of (2.18) with*

$$u = k(x, z) \tag{2.21}$$

$$\dot{z} = \Gamma e^{-z}\left(V(x,z) - \frac{\varepsilon^2}{2}\right)^+ \quad (2.22)$$

and $(x(0), z(0)) = (x_0, z_0)$ that corresponds to $d \in L^\infty(\mathbb{R}_+; \mathbb{R}^l)$, $\theta \in L^\infty(\mathbb{R}_+; \Theta)$ exists for all $t \geq 0$ and satisfies the following estimates for all $t \geq 0$:

$$\limsup_{t \to +\infty}(|x_1(t)|) \leq \varepsilon, \quad \limsup_{t \to +\infty}(V(x(t), z(t))) \leq \frac{\varepsilon^2}{2} \quad (2.23)$$

$$z_0 \leq z(t) \leq \lim_{s \to +\infty}(z(s)) \leq B(x_0, z_0, \|d\|_\infty, \|\theta\|_\infty) \quad (2.24)$$

$$|x(t)| \leq B(x_0, z_0, \|d\|_\infty, \|\theta\|_\infty) \quad (2.25)$$

$$|x_1(t)|^2 \leq 2V(x(t), z(t)) \leq 2e^{-ct}V(x_0, z_0) + 2a\frac{\|d\|_\infty^2 + \left(\left(\|\theta\|_\infty - b - \lambda(e^{z_0})\right)^+\right)^2}{c(1 + \kappa(e^{z_0}))} \quad (2.26)$$

Moreover, if $\lim_{t \to +\infty}(d(t)) = 0$ and $\limsup_{t \to +\infty}(|\theta(t)|) \leq b + \lambda\left(e^{\lim_{t \to +\infty}(z(t))}\right)$ then $\lim_{t \to +\infty}(|x(t)|) = 0$.

## 3. Illustrative Example

The following system is inspired by the aircraft wing rock example on pages 180-183 in the book [15]:

$$\begin{aligned}\dot{x}_1 &= x_2 \\ \dot{x}_2 &= \theta_1 x_1 + \theta_2 x_2 + \theta_3 x_1 x_2 + \theta_4 x_2^2 + x_3 + d_1 \\ \dot{x}_3 &= u + d_2\end{aligned} \quad (3.1)$$

where $x = (x_1, x_2, x_3) \in \mathbb{R}^3$ is the state, $u \in \mathbb{R}$ is the control input, $d = (d_1, d_2) \in \mathbb{R}^2$, $\theta = (\theta_1, \theta_2, \theta_3, \theta_4) \in \mathbb{R}^4$ are disturbances. Applying Lemma 4 in Section 4 repeatedly, we can guarantee that for every $\Gamma > 0, c \geq 1/2, K \geq 28c, \varepsilon > 0$ the feedback law (DADS controller)

$$\begin{aligned}u &= -\left(2c + K\rho^2\left(L(x) + 4\zeta x_2^3\right)\right)x_3 - \zeta - x_2 - 2\Gamma K\rho L(x)(V(x,z) - r)^+ \zeta \\ &\quad -K\rho^2 x_2\left(4x_1^3 \zeta + 2cL(x)\right) - 42c(2c+1)\rho^2 L(x)\left(1 + 18cK\rho^2 L(x)\right)^2 \xi \\ \dot{z} &= \Gamma e^{-z}(V(x,z) - \varepsilon)^+\end{aligned} \quad (3.2)$$

where

$$\zeta := x_2 + 2cx_1$$
$$\rho := 1 + e^z$$
$$L(x) := 1 + x_1^4 + x_2^4 \quad (3.3)$$
$$\xi := x_3 + x_1 + 2cx_2 + K\rho^2 L(x)\zeta$$
$$V(x,z) := \frac{1}{2}x_1^2 + \frac{1}{2}\zeta^2 + \frac{1}{2}\xi^2$$

satisfies the following inequality for all $(x,z) \in \mathbb{R}^3 \times \mathbb{R}, \theta \in \mathbb{R}^4, d \in \mathbb{R}^2$:

$$\dot{V} \leq -cV(x,z) + 2\frac{|d|^2 + \left(\left(|\theta| - 1 - e^z\right)^+\right)^2}{1+e^z} \quad (3.4)$$

Using definitions (3.3) and estimate (3.4) we can apply Theorem 2 with $b=1$, $a=2$, $\kappa(s) = \lambda(s) = s$ for $s \geq 0$ and guarantee estimates (2.13)-(2.17) for an appropriate function $B \in C^0(\mathbb{R}^3 \times \mathbb{R} \times \mathbb{R}_+ \times \mathbb{R}_+)$. Using estimates (2.13), (2.16) in conjunction with definitions (3.3) (which guarantees that $V(x,z) \geq \frac{\sqrt{c^2+1}-c}{2(\sqrt{c^2+1}+c)}|(x_1,x_2)|^2$ for all $(x,z) \in \mathbb{R}^3 \times \mathbb{R}$), we obtain the following attractivity property for all $\theta, d \in L^\infty(\mathbb{R}_+)$ and for all initial conditions

$$\limsup_{t \to +\infty}\left(|(x_1(t), x_2(t))|\right) \leq \left(\sqrt{c^2+1}+c\right)\sqrt{2\varepsilon}$$

as well as the following IOS-like property with $Y = (x_1, x_2) \in \mathbb{R}^2$ as output

$$|(x_1(t), x_2(t))|^2 \leq 2\left(\sqrt{c^2+1}+c\right)^2 e^{-ct} V(x(0), z(0))$$
$$+ \frac{4\left(\sqrt{c^2+1}+c\right)^2}{c\left(1+e^{z(0)}\right)}\left(\|d\|_\infty^2 + \left(\left(\|\theta\|_\infty - 1 - e^{z(0)}\right)^+\right)^2\right)$$

For comparison purposes we give the adaptive controller that is designed for the case of constant $\theta$ using the methodology proposed in Chapter 4 of [15] with $\sigma-$ modification:

$$u = -(w_1 + 2c + w_3 x_2)x_1 - (\hat{\theta}_1 + w_2 + 2 + 2Kc)x_2 - (\hat{\theta}_3 + w_4)x_2^2$$
$$-\varphi(\hat{\theta}_1 x_1 + \hat{\theta}_2 x_2 + \hat{\theta}_3 x_1 x_2 + \hat{\theta}_4 x_2^2 + x_3) - (K + \varphi^2)\chi \quad (3.5)$$

$$\frac{d\hat{\theta}_i}{dt} = w_i, \quad i = 1, 2, 3, 4$$

$$w_1 = \Gamma(\zeta + \varphi\chi)x_1 - \sigma\hat{\theta}_1$$

$$w_2 = \Gamma(\zeta + \varphi\chi)x_2 - \sigma\hat{\theta}_2 \quad (3.6)$$

$$w_3 = \Gamma(\zeta + \varphi\chi)x_1 x_2 - \sigma\hat{\theta}_3$$

$$w_4 = \Gamma(\zeta + \varphi\chi)x_2^2 - \sigma\hat{\theta}_4$$

where $c, \Gamma > 0$, $K \geq 1 + 2c$, $\sigma \geq 0$ are constants and

$$\zeta := x_2 + 2cx_1$$

$$\chi := \hat{\theta}_1 x_1 + \hat{\theta}_2 x_2 + \hat{\theta}_3 x_1 x_2 + \hat{\theta}_4 x_2^2 + 2cx_2 + K\zeta + x_1 + x_3 \quad (3.7)$$

$$\varphi := 2c + K + \hat{\theta}_2 + \hat{\theta}_3 x_1 + 2\hat{\theta}_4 x_2$$

The adaptive feedback law (3.5), (3.6) guarantees the following inequality for the function $W(x,\hat{\theta}) = \frac{1}{2}x_1^2 + \frac{1}{2}\zeta^2 + \frac{1}{2}\chi^2 + \frac{1}{2\Gamma}\sum_{i=1}^{4}(\hat{\theta}_i - \theta_i)^2$ and for all $x \in \mathbb{R}^3$, $d \in \mathbb{R}^2$, $\hat{\theta} \in \mathbb{R}^4$, $\theta \in \mathbb{R}^4$:

$$\dot{W} \leq -c(x_1^2 + \zeta^2 + \chi^2) - \frac{\sigma}{2\Gamma}\sum_{i=1}^{4}(\hat{\theta}_i - \theta_i)^2 + \frac{1}{2}|d|^2 + \frac{\sigma}{2\Gamma}|\theta|^2 \quad (3.8)$$

Next, we present plots for the following parameters and initial conditions:

$$c = 0.5, K = 14, \Gamma = 20, \varepsilon = 0.01, \theta_1 = \theta_2 = 20, \theta_3 = 2, \theta_4 = 1,$$

$$x_1(0) = 1, \ x_2(0) = -0.5, \ x_3(0) = -18, \ z(0) = -\ln(10), \ \hat{\theta}(0) = 0 \in \mathbb{R}^4$$

Figure 5 shows the problem with the use of the adaptive feedback law (3.5), (3.6) with no $\sigma$-modification (i.e., with $\sigma = 0$): the controller state $\hat{\theta} \in \mathbb{R}^4$ grows without bound when a persistent disturbance is present. Otherwise the adaptive feedback law (3.5), (3.6) with no $\sigma$-modification achieves a very good performance: see the output regulation properties for the output $Y = (x_1, x_2) \in \mathbb{R}^2$ in Figure 1 (no disturbance) and in Figure 4 (case of persistent disturbance). The reader can also see in Figure 2 that the norm of the controller state $|\hat{\theta}|$ remains bounded when no disturbance is present. However, as noticed above this is not true when a persistent disturbance is present and therefore we are forced to use $\sigma$-modification.

The comparison between the performance induced by the DADS controller (3.2) and adaptive feedback law (3.5), (3.6) with $\sigma$-modification is made in Figure 1, Figure 2, Figure 3 for the disturbance-free case and Figure 4, Figure 5, Figure 6 for the case of a persistent disturbance. The DADS controller does not only achieve better regulation of the output $Y = (x_1, x_2) \in \mathbb{R}^2$ but also achieves a smaller control effort (both in the control input and the controller states). This is true both for the disturbance-free case and the case of a persistent disturbance.

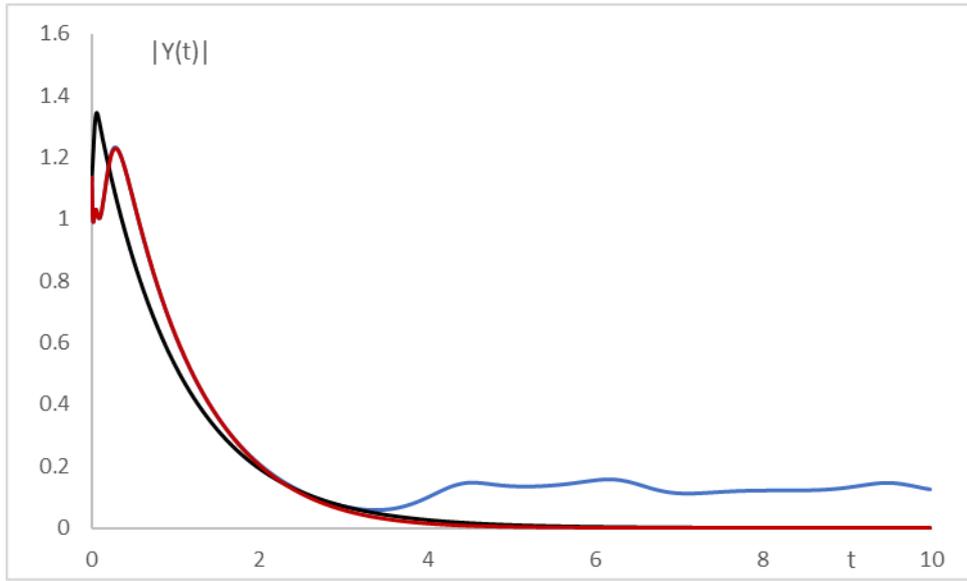

**Fig. 1:** Plot of $|Y(t)| = \sqrt{x_1^2(t) + x_2^2(t)}$ for the disturbance-free case ($d_1(t) = d_2(t) \equiv 0$). Black line for system (3.1) with DADS feedback (3.2); blue line for system (3.1) with $\sigma$-adaptive feedback (3.5), (3.6) and $\sigma = 0.4$; red line for system (3.1) with (3.5), (3.6) and $\sigma = 0$.

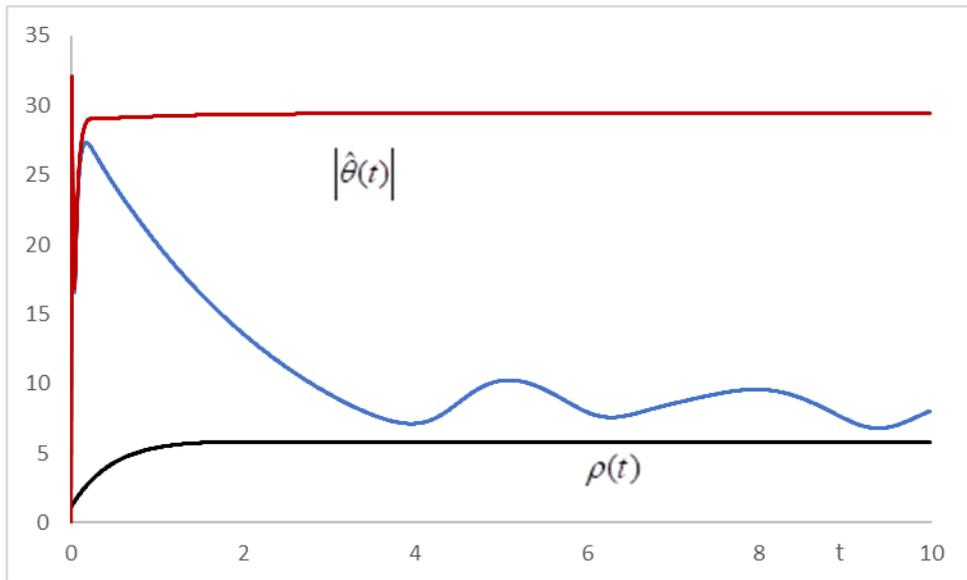

**Fig. 2:** Plot of $\rho(t) = 1 + e^{z(t)}$ for system (3.1) with DADS feedback (3.2) (black line), plot of $|\hat{\theta}(t)|$ for system (3.1) with $\sigma$-adaptive feedback (3.5), (3.6) and $\sigma = 0.4$ (blue line) and plot of $|\hat{\theta}(t)|$ for system (3.1) with (3.5), (3.6) and $\sigma = 0$ (red line). All plots are for the disturbance-free case (i.e., for the case $d_1(t) = d_2(t) \equiv 0$).

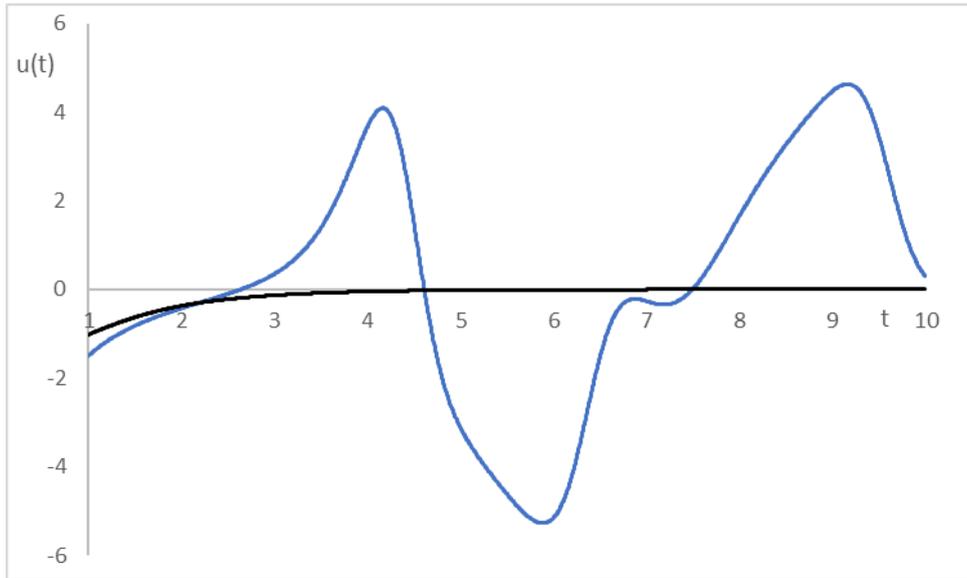

**Fig. 3:** Plot of $u(t)$ for $t \in [1,10]$ and the disturbance-free case ($d_1(t) = d_2(t) \equiv 0$). Black line for system (3.1) with DADS feedback (3.2); blue line for system (3.1) with $\sigma-$ adaptive feedback (3.5), (3.6) and $\sigma = 0.4$.

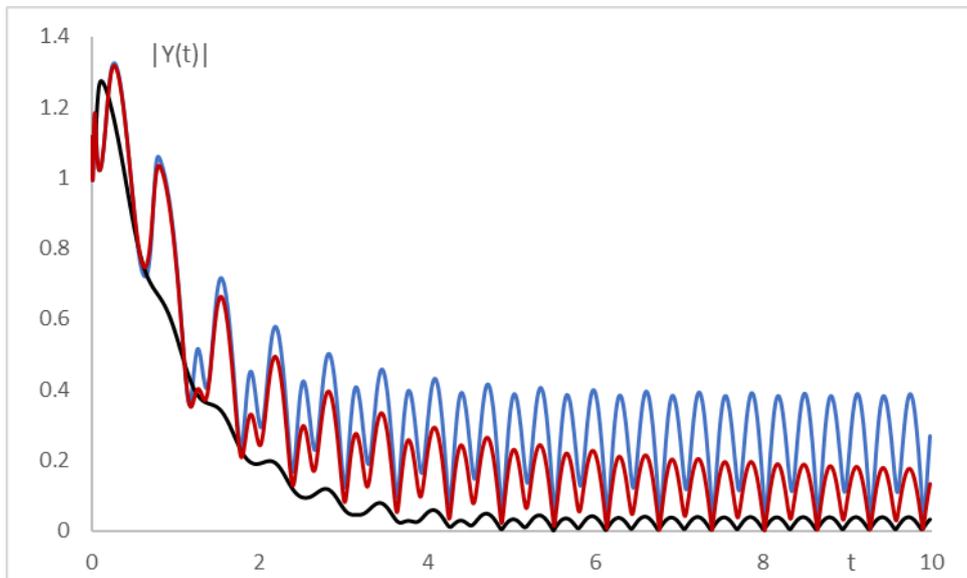

**Fig. 4:** Plot of $|Y(t)| = \sqrt{x_1^2(t) + x_2^2(t)}$ for the case $d_1(t) = 20\cos(10t)$, $d_2(t) = 10\cos(20t)$. Black line for system (3.1) with DADS feedback (3.2); blue line for system (3.1) with $\sigma-$ adaptive feedback (3.5), (3.6) and $\sigma = 0.4$; red line for system (3.1) with (3.5), (3.6) and $\sigma = 0$.

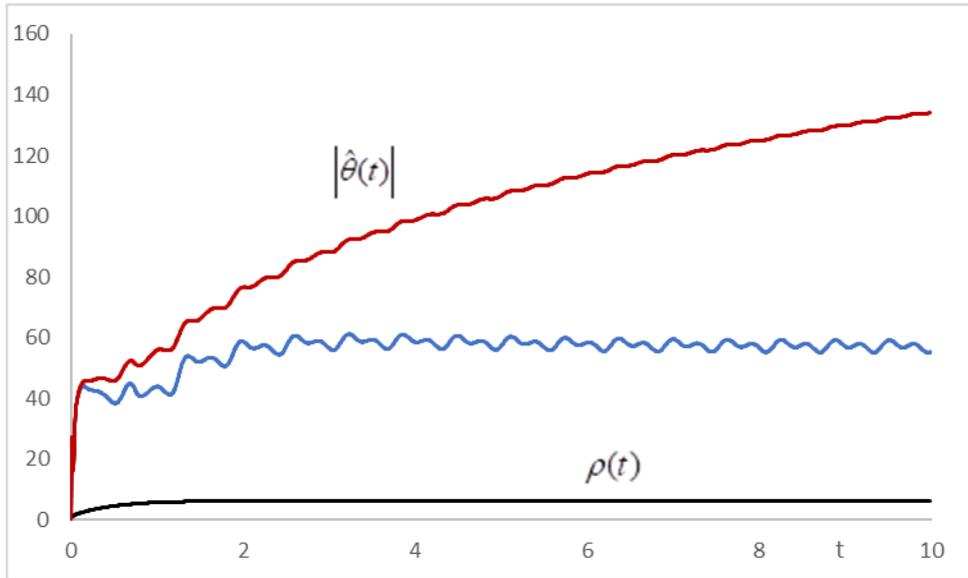

**Fig. 5:** Plot of $\rho(t) = 1 + e^{z(t)}$ for system (3.1) with DADS feedback (3.2) (black line), plot of $|\hat{\theta}(t)|$ for system (3.1) with (3.5), (3.6) and $\sigma = 0.4$ (blue line) and plot of drifting in $|\hat{\theta}(t)|$ for system (3.1) with $\sigma-$adaptive feedback (3.5), (3.6) and $\sigma = 0$ (red line). All plots are for the case $d_1(t) = 20\cos(10t)$, $d_2(t) = 10\cos(20t)$.

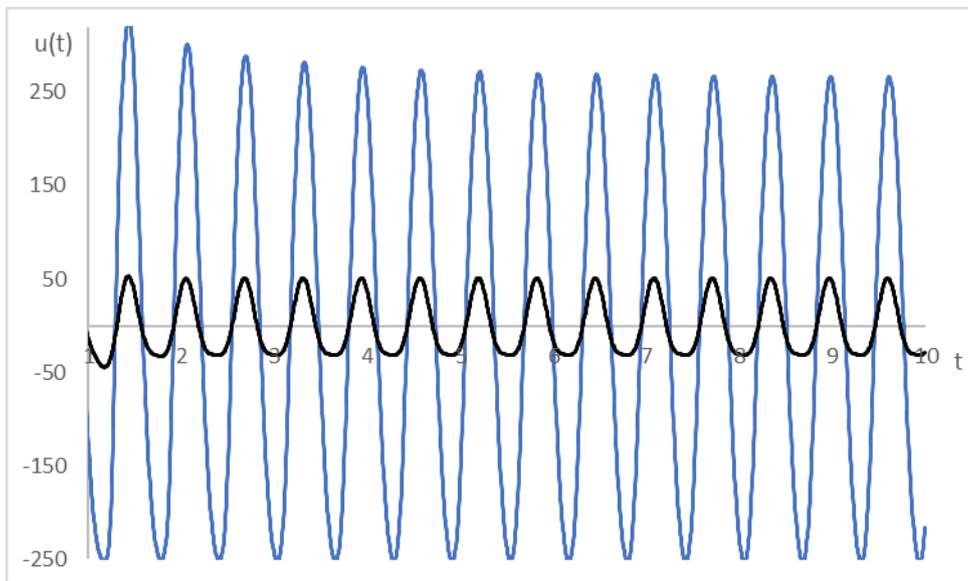

**Fig. 6:** Plot of the input $u(t)$ for $t \in [1,10]$ and the case $d_1(t) = 20\cos(10t)$, $d_2(t) = 10\cos(20t)$. Black line for system (3.1) with DADS feedback (3.2); blue line for system (3.1) with $\sigma-$adaptive feedback (3.5), (3.6) and $\sigma = 0.4$.

## 4. Proofs

We first provide the proof of Theorem 2.

**Proof of Theorem 2:** Let arbitrary $(x_0, z_0) \in \mathbb{R}^n \times \mathbb{R}$, $d \in L^\infty(\mathbb{R}_+; D)$, $\theta \in L^\infty(\mathbb{R}_+; \Theta)$ be given. Since the functions $f : \mathbb{R}^n \times \mathbb{R} \times \Theta \times D \to \mathbb{R}^n$ and $V$ are locally Lipschitz there exists an unique solution $(x(t), z(t))$ of (2.11) with $(x(0), z(0)) = (x_0, z_0)$ defined on a maximal interval $[0, t_{max})$, where $t_{max} \in (0, +\infty]$. For every $t \in [0, t_{max})$, we have from (2.11) that $\dot{z}(t) \geq 0$. Therefore, $z(t)$ is non-decreasing on $[0, t_{max})$. Moreover, the mapping

$$z \to \frac{|d|^2 + \left(\left(|\theta| - b - \lambda(e^z)\right)^+\right)^2}{1 + \kappa(e^z)}$$

is non-increasing. Using (2.12) we conclude that the following inequalities hold almost everywhere in $[0, t_{max})$:

$$\frac{d}{dt} V(x(t), z(t)) \leq -cV(x(t), z(t)) + a \frac{|d(t)|^2 + \left(\left(|\theta(t)| - b - \lambda\left(e^{z(t)}\right)\right)^+\right)^2}{1 + \kappa(e^{z(t)})} \quad (4.1)$$

$$\frac{d}{dt} V(x(t), z(t)) \leq -cV(x(t), z(t)) + a \frac{\|d\|_\infty^2 + \left(\left(\|\theta\|_\infty - b - \lambda\left(e^{z_0}\right)\right)^+\right)^2}{1 + \kappa(e^{z_0})} \quad (4.2)$$

Integrating the differential inequality (4.2), we obtain for all $t \in [0, t_{max})$:

$$V(x(t), z(t)) \leq e^{-ct} V(x_0, z_0) + a \frac{\|d\|_\infty^2 + \left(\left(\|\theta\|_\infty - b - \lambda\left(e^{z_0}\right)\right)^+\right)^2}{c\left(1 + \kappa(e^{z_0})\right)} \quad (4.3)$$

Next, we show that $z(t)$ is bounded. Since $z(t)$ is non-decreasing on $[0, t_{max})$, we have $z(t) \geq z_0$ for all $t \in [0, t_{max})$. Define

$$s := V(x_0, z_0) + |z_0| + \|\theta\|_\infty^2 + \|d\|_\infty^2 \quad (4.4)$$

$$A(r) := \kappa^{-1}\left(\left(\frac{ar}{c\varepsilon} - 1\right)^+\right), \text{ for } r \geq 0 \quad (4.5)$$

$$G(r) := \ln\left(1 + \max\{A(r), e^r\} + \frac{\Gamma}{c}\left(1 + \frac{a}{c}\right)r\right), \text{ for } r \geq 0 \quad (4.6)$$

We next distinguish the following cases:

*Case 1:* $e^{z(t)} \leq 1 + \max\{A(s), e^s\}$ for all $t \in [0, t_{max})$

In this case, it follows from (4.6) that:

$$z(t) \leq G(s), \text{ for all } t \in [0, t_{max})$$

*Case 2:* There exists $T \in (0, t_{max})$ such that $e^{z(T)} > 1 + \max\{A(s), e^s\}$. Since $e^{z(0)} < e^s < 1 + \max\{A(s), e^s\}$, continuity of $z(t)$ guarantees that there exist $t_0 \in (0, T)$ such that:

$$e^{z(t_0)} = 1 + \max\{A(s), e^s\} \tag{4.7}$$

Since $z(t)$ is non-decreasing, it follows that:

$$e^{z(t)} \geq 1 + \max\{A(s), e^s\}, \text{ for } t \in [t_0, t_{max}) \tag{4.8}$$

Inequality (4.8) in conjunction with (4.4) and (4.5) gives the following inequalities for $t \in [t_0, t_{max})$:

$$e^{z(t)} \geq A(s) \geq \kappa^{-1}\left(\left(\frac{as}{c\varepsilon} - 1\right)^+\right) \geq \kappa^{-1}\left(\left(a\frac{\|\theta\|_\infty^2 + \|d\|_\infty^2}{c\varepsilon} - 1\right)^+\right)$$

$$\Rightarrow \kappa\left(e^{z(t)}\right) \geq \left(a\frac{\|\theta\|_\infty^2 + \|d\|_\infty^2}{c\varepsilon} - 1\right)^+ \geq a\frac{\|\theta\|_\infty^2 + \|d\|_\infty^2}{c\varepsilon} - 1 \Rightarrow c\varepsilon \geq a\frac{\|\theta\|_\infty^2 + \|d\|_\infty^2}{1 + \kappa\left(e^{z(t)}\right)}$$

Consequently, we obtain from (4.2) for almost all $t \in [t_0, t_{max})$:

$$\frac{d}{dt}V(x(t), z(t)) \leq -cV(x(t), z(t)) + c\varepsilon$$

Integrating the above differential inequality, we obtain for all $t \in [t_0, t_{max})$:

$$V(x(t), z(t)) \leq e^{-c(t-t_0)}V(x(t_0), z(t_0)) + \varepsilon \tag{4.9}$$

Moreover, from estimate (4.3) and definition (4.4), we get the following estimate for all $t \in [0, t_{max})$:

$$V(x(t), z(t)) \leq V(x_0, z_0) + \frac{a}{c}\left(\|\theta\|_\infty^2 + \|d\|_\infty^2\right) \leq s + \frac{a}{c}s \tag{4.10}$$

Using (2.11) and inequality (4.9), we have for all $t \in [t_0, t_{max})$:

$$\frac{d}{dt}\left(e^{z(t)}\right) \leq \Gamma e^{-c(t-t_0)}V(x(t_0), z(t_0)) \tag{4.11}$$

Integrating (4.11) and using (4.10) and (4.7), we obtain for all $t \in [t_0, t_{max})$:

$$e^{z(t)} \leq e^{z(t_0)} + \frac{\Gamma}{c} V(x(t_0), z(t_0))$$

$$= 1 + \max\{A(s), e^s\} + \frac{\Gamma}{c} V(x(t_0), z(t_0)) \quad (4.12)$$

$$\leq 1 + \max\{A(s), e^s\} + \frac{\Gamma}{c}\left(1 + \frac{a}{c}\right) s$$

Since $z(t)$ is non-decreasing, inequality (4.12) holds for all $t \in [0, t_{max})$.

It follows from (4.12) and definitions (4.4), (4.6) that in both cases the following estimates hold for all $t \in [0, t_{max})$:

$$z_0 \leq z(t) \leq G(s) \quad (4.13)$$

$$|z(t)| \leq G(s) \quad (4.14)$$

Estimate (4.14) implies that $z(t)$ is bounded. Therefore, assumption (i) guarantees that $x(t)$ is bounded for all $t \in [0, t_{max})$. Consequently, $t_{max} = +\infty$. Moreover, estimate (4.3) implies directly the required estimate (2.16).

Since $z(t)$ is non-decreasing and bounded from above, $\lim_{\tau \to +\infty} (z(\tau))$ exists. Thus, we obtain from (4.13):

$$z_0 \leq z(t) \leq \lim_{\tau \to +\infty} (z(\tau)) \leq G(s) \quad (4.15)$$

Next, define the following function for $r \geq 0$:

$$H(r) := \max\left\{|y| : V(y, \xi) \leq \left(1 + \frac{a}{c}\right) r \text{ and } |\xi| \leq G(r)\right\} \quad (4.16)$$

Assumption (i) guarantees that $H$ is well-defined. Since $G$ is non-decreasing (recall definition (4.6)), it follows from (A16) that $H$ is non-decreasing. Moreover, it follows from estimates (4.10), (4.14) and definition (4.16):

$$|x(t)| \leq H(s), \text{ for all } t \geq 0 \quad (4.17)$$

Define the continuous function $\tilde{H}(r) := G(r) + \int_r^{r+1} H(l) dl$ for $r \geq 0$ (notice that since $H$ is non-decreasing, it is locally bounded and Riemann-integrable on every bounded interval of $\mathbb{R}_+$). We conclude from (4.17), (4.13) and (4.4) that estimates (2.14), (2.15) hold with $B(x_0, z_0, \|d\|_\infty, \|\theta\|_\infty) = \tilde{H}(V(x_0, z_0) + |z_0| + \|\theta\|_\infty^2 + \|d\|_\infty^2)$.

We next show estimate (2.13). Estimate (4.13) and the fact that $z(t)$ is non-decreasing guarantees that the function $z(t)$ has a finite limit as $t \to +\infty$. This implies that the function $e^{z(t)}$ has a finite limit as $t \to +\infty$. Moreover, the fact that $d \in L^\infty(\mathbb{R}_+; D)$, $\theta \in L^\infty(\mathbb{R}_+; \Theta)$ and (4.14), (4.17) implies that $\dot{V}(t)$ is of class $L^\infty(\mathbb{R}_+)$. It follows that the function $\frac{d}{dt}(e^{z(t)}) = \Gamma(V(x(t), z(t)) - \varepsilon)^+$ is uniformly continuous. From Barbalat's Lemma (see [9]), we have:

$$\lim_{t \to +\infty} \left( \frac{d}{dt}(e^{z(t)}) \right) = \lim_{t \to +\infty} \left( \Gamma(V(x(t), z(t)) - \varepsilon)^+ \right) = 0 \tag{4.18}$$

Therefore, estimate (2.13) holds.

Next, we show estimate (2.17). Let arbitrary $\epsilon > 0$ be given. Then there exists $\bar{T} > 0$ such that

$$|d(t)| \leq \limsup_{\tau \to +\infty}(|d(\tau)|) + \epsilon \text{ for } t \geq \bar{T} \text{ a.e.}$$

$$|\theta(t)| \leq \limsup_{\tau \to +\infty}(|\theta(\tau)|) + \epsilon \text{ for } t \geq \bar{T} \text{ a.e.}$$

$$z(t) \geq \lim_{\tau \to +\infty}(z(\tau)) - \epsilon \text{ for all } t \geq \bar{T}$$

Define $L_d = \limsup_{\tau \to +\infty}(|d(\tau)|)$, $L_\theta = \limsup_{\tau \to +\infty}(|\theta(\tau)|)$ and $L_z = \lim_{\tau \to +\infty}(z(\tau))$. It follows from (4.1) that the following differential inequality holds for $t \geq \bar{T}$ a.e.:

$$\frac{d}{dt}V(x(t), z(t)) \leq -cV(x(t), z(t)) + a \frac{(L_d + \epsilon)^2 + \left(\left(L_\theta + \epsilon - b - \lambda\left(e^{L_z - \epsilon}\right)\right)^+\right)^2}{1 + \kappa\left(e^{L_z - \epsilon}\right)} \tag{4.19}$$

Integrating (4.19), we get for all $t \geq \bar{T}$:

$$V(x(t), z(t)) \leq e^{-c(t - \bar{T})} V(x(\bar{T}), z(\bar{T})) + a \frac{(L_d + \epsilon)^2 + \left(\left(L_\theta + \epsilon - b - \lambda\left(e^{L_z - \epsilon}\right)\right)^+\right)^2}{c\left(1 + \kappa\left(e^{L_z - \epsilon}\right)\right)} \tag{4.20}$$

Consequently, we get from estimate (A20):

$$\limsup_{t \to +\infty}(V(x(t), z(t))) \leq a \frac{(L_d + \epsilon)^2 + \left(\left(L_\theta + \epsilon - b - \lambda\left(e^{L_z - \epsilon}\right)\right)^+\right)^2}{c\left(1 + \kappa\left(e^{L_z - \epsilon}\right)\right)} \tag{4.21}$$

Since (4.21) holds for arbitrary $\epsilon > 0$, we obtain the required estimate (2.17).

The proof is complete. ◁

The proofs of Theorem 1 and Theorem 2 are based on the following lemma which allows a step-by-step construction of the DADS controller.

**Lemma 4 (A Backstepping Lemma):** *Consider the control system*

$$\dot{x} = f(x) + g(x,\theta)y + \Phi(x)\theta + G(x)d$$
$$\dot{y} = h(x,y) + \tilde{g}(x,y,\theta)u + \varphi'(x,y)\theta + \alpha'(x,y)d \tag{4.22}$$

*where* $(x,y) \in \mathbb{R}^n \times \mathbb{R}$ *is the state,* $u \in \mathbb{R}$ *is the control input,* $d \in \mathbb{R}^l$, $\theta \in \Theta \subseteq \mathbb{R}^p$ *are disturbances,* $f : \mathbb{R}^n \to \mathbb{R}^n$, $g : \mathbb{R}^n \times \Theta \to \mathbb{R}^n$, $\Phi : \mathbb{R}^n \to \mathbb{R}^{n \times p}$, $G : \mathbb{R}^n \to \mathbb{R}^{n \times l}$, $h : \mathbb{R}^n \times \mathbb{R} \to \mathbb{R}$, $\varphi : \mathbb{R}^n \times \mathbb{R} \to \mathbb{R}^p$, $\alpha : \mathbb{R}^n \times \mathbb{R} \to \mathbb{R}^l$, $\tilde{g} : \mathbb{R}^n \times \mathbb{R} \times \Theta \to \mathbb{R}$ *are smooth mappings with* $f(0) = 0$, $\Phi(0) = 0$, $h(0,0) = 0$ *and* $\varphi(0,0) = 0$.

*Assume that there exist constants* $a, b, c, \Gamma, \varepsilon, \eta > 0$, *smooth functions* $\kappa, \lambda \in K_\infty$, $\sigma, k, V \in C^\infty(\mathbb{R}^n \times \mathbb{R})$, $\eta \in C^\infty(\mathbb{R}^n \times \mathbb{R}; (0, +\infty))$, $\mu \in C^\infty(\mathbb{R}^n)$ *with* $V(0,z) = k(0,z) = 0$ *for* $z \in \mathbb{R}$ *and* $V(x,z) > 0$ *for* $x \neq 0$, $z \in \mathbb{R}$, *such that the following inequalities hold for all* $\theta \in \Theta$, $d \in \mathbb{R}^l$, $(x,y,z) \in \mathbb{R}^n \times \mathbb{R} \times \mathbb{R}$:

$$|x|^2 \leq \sigma(x,z)V(x,z) \tag{4.23}$$

$$\frac{\partial V}{\partial x}(x,z)\big(f(x) + g(x,\theta)k(x,z) + \Phi(x)\theta + G(x)d\big)$$
$$+ \frac{\partial V}{\partial z}(x,z)\Gamma e^{-z}\big(V(x,z) - \varepsilon\big)^+ \leq -cV(x,z) + a\frac{|d|^2 + \big((|\theta| - b - \lambda(e^z))^+\big)^2}{1 + \kappa(e^z)} \tag{4.24}$$

$$\tilde{g}(x,y,\theta) \geq \eta(x,y) \tag{4.25}$$

$$|g(x,\theta)| \leq \mu(x)(1 + |\theta|) \tag{4.26}$$

*Define the following function for all* $(x,y,z) \in \mathbb{R}^n \times \mathbb{R} \times \mathbb{R}$:

$$\bar{V}(x,y,z) = V(x,z) + \frac{1}{2}(y - k(x,z))^2 \tag{4.27}$$

*Then there exist functions* $\bar{\sigma}, \bar{k} \in C^\infty(\mathbb{R}^n \times \mathbb{R} \times \mathbb{R})$ *with* $\bar{k}(0,0,z) = 0$ *for* $z \in \mathbb{R}$, *such that the following inequalities hold for all* $d \in \mathbb{R}^l$, $\theta \in \Theta$ *and* $(x,y,z) \in \mathbb{R}^n \times \mathbb{R} \times \mathbb{R}$:

$$|x|^2 + |y|^2 \leq \bar{\sigma}(x,y,z)\bar{V}(x,y,z) \tag{4.28}$$

$$\frac{\partial \bar{V}}{\partial x}(x,y,z)\big(f(x)+g(x,\theta)y+\Phi(x)\theta+G(x)d\big)+\frac{\partial \bar{V}}{\partial z}(x,y,z)\Gamma e^{-z}\big(\bar{V}(x,y,z)-\varepsilon\big)^{+}$$

$$+\frac{\partial \bar{V}}{\partial y}(x,y,z)\big(h(x,y)+\tilde{g}(x,y,\theta)\bar{k}(x,y,z)+\varphi'(x,y)\theta+\alpha'(x,y)d\big) \quad (4.29)$$

$$\leq -\frac{c}{2}\bar{V}(x,y,z)+2a\frac{|d|^{2}+\big((|\theta|-b-\lambda(e^{z}))^{+}\big)^{2}}{1+\kappa(e^{z})}$$

**Remark:** The function $\bar{V} \in C^{\infty}(\mathbb{R}^{n}\times\mathbb{R}\times\mathbb{R})$ defined by (4.27) satisfies $\bar{V}(0,0,z)=0$ for $z\in\mathbb{R}$ and $\bar{V}(x,y,z)>0$ for $(x,y)\neq 0$, $z\in\mathbb{R}$. Moreover, if the set $\{(x,z)\in\mathbb{R}^{n}\times[-M,M]:V(x,z)\leq M\}$ is bounded for each $M>0$ then the set $\{(x,y,z)\in\mathbb{R}^{n}\times\mathbb{R}\times[-M,M]:\bar{V}(x,y,z)\leq M\}$ is bounded for each $M>0$.

**Proof of Lemma 4:** Without loss of generality we may assume that $\sigma(x,z)>0$ for all $x\in\mathbb{R}^{n}$, $z\in\mathbb{R}$.

We start by using the fact that for every smooth function $\psi:\mathbb{R}^{n}\times\mathbb{R}\to\mathbb{R}^{m}$ with $\psi(0,z)=0$ for all $z\in\mathbb{R}$ there exists a smooth positive function $L_{\psi}:\mathbb{R}^{n}\times\mathbb{R}\to(0,+\infty)$ such that the inequality $|\psi(x,z)|\leq L_{\psi}(x,z)|x|$ holds for all $(x,z)\in\mathbb{R}^{n}\times\mathbb{R}$.

Since $V(0,z)=0$ for all $z\in\mathbb{R}$ and $V(x,z)>0$ for all $x\neq 0$, $z\in\mathbb{R}$, it follows that for each fixed $z\in\mathbb{R}$, the mapping $x\to V(x,z)$ has a global minimum at $x=0$. Consequently, $\frac{\partial V}{\partial x}(0,z)=0$ for all $z\in\mathbb{R}$. Moreover, since $k(0,z)=0$ for all $z\in\mathbb{R}$ and since $f(0)=0$, $\Phi(0)=0$, $h(0,0)=0$, $\varphi(0,0)=0$, we conclude that there exist smooth positive functions $R:\mathbb{R}^{n}\times\mathbb{R}\to(0,+\infty)$, $r:\mathbb{R}^{n}\to(0,+\infty)$, $\rho:\mathbb{R}^{n}\times\mathbb{R}\to(0,+\infty)$ such that the following inequalities hold for all $(x,z)\in\mathbb{R}^{n}\times\mathbb{R}$ and $y\in\mathbb{R}$:

$$\left|\frac{\partial V}{\partial x}(x,z)\right|+|k(x,z)|+\left|\frac{\partial k}{\partial x}(x,z)\Phi(x)\right|\leq R(x,z)|x|$$
$$|f(x)|\leq r(x)|x| \quad (4.30)$$
$$|h(x,y)|+|\varphi(x,y)|\leq \rho(x,y)|x|+\rho(x,y)|y|$$

Define for all $(x,z)\in\mathbb{R}^{n}\times\mathbb{R}$ and $y\in\mathbb{R}$

$$\bar{k}(x,y,z):=-\frac{M(x,y,z)}{\eta(x,y)}(y-k(x,z)) \quad (4.31)$$

$$\bar{\sigma}(x,y,z):=\big(1+2R^{2}(x,z)\big)\sigma(x,z)+4 \quad (4.32)$$

where $\eta$ is the function appearing in (4.25) and $M : \mathbb{R}^n \times \mathbb{R} \times \mathbb{R} \to (0, +\infty)$ is defined by the following formulas:

$$M(x, y, z) := \frac{c}{4} + \frac{\Gamma^2 e^{-2z}}{4c}\left(1 + \left(\frac{\partial k}{\partial z}(x, z)\right)^2\right)^2 V(x, z) + P(x, y, z)\left(b + \lambda\left(e^z\right)\right)$$

$$+ \frac{1}{2}\left(\frac{\Gamma e^{-z}}{4}\left(1 + (y - k(x, z))^2\right) + \mu(x)\right)\left(1 + \left(\frac{\partial k}{\partial z}(x, z)\right)^2\right) + \rho(x, y)$$

$$+ \frac{\sigma(x, z)}{c} P^2(x, y, z)\left(b + 1 + \lambda\left(e^z\right)\right)^2 + \frac{1 + \kappa\left(e^z\right)}{4a}\left|\alpha'(x, y) - \frac{\partial k}{\partial x}(x, z) G(x)\right|^2$$

$$+ \frac{1 + \kappa\left(e^z\right)}{2a} P^2(x, y, z)\left((y - k(x, z))^2 + |x|^2\right) + \frac{\Gamma e^{-z}}{4}\left(1 + \left(\frac{\partial V}{\partial z}(x, z)\right)^2\right) \quad (4.33)$$

$$P(x, y, z) := \frac{r(x) + \mu(x)}{2}\left(1 + \left|\frac{\partial k}{\partial x}(x, z)\right|^2\right) + \rho(x, y)$$

$$+ \left(1 + \frac{\mu(x)}{2}\left(3 + \left|\frac{\partial k}{\partial x}(x, z)\right|^2\right) + \rho(x, y)\right) R(x, z) \quad (4.34)$$

In what follows, we denote by $\Lambda$ the quantity appearing in the left hand side of (4.29) and $d \in \mathbb{R}^l$, $\theta \in \Theta$, $(x, y, z) \in \mathbb{R}^n \times \mathbb{R} \times \mathbb{R}$ are arbitrary given vectors. We also use the notation $s = y - k(x, z)$. Using definition (4.27) we get:

$$\Lambda = \left(\frac{\partial V}{\partial z}(x, z) - s \frac{\partial k}{\partial z}(x, z)\right) \Gamma e^{-z} \left(\overline{V}(x, y, z) - \varepsilon\right)^+$$

$$+ \left(\frac{\partial V}{\partial x}(x, z) - s \frac{\partial k}{\partial x}(x, z)\right)\left(f(x) + g(x, \theta) y + \Phi(x)\theta + G(x)d\right) \quad (4.35)$$

$$+ s\left(h(x, y) + \tilde{g}(x, y, \theta)\overline{k}(x, y, z) + \varphi'(x, y)\theta + \alpha'(x, y)d\right)$$

Inequality (4.24) in conjunction with (4.35) gives:

$$\Lambda \leq -cV(x,z) - s\frac{\partial k}{\partial z}(x,z)\Gamma e^{-z}\left(\overline{V}(x,y,z) - \varepsilon\right)^{+}$$

$$+ \frac{\partial V}{\partial z}(x,z)\Gamma e^{-z}\left(\left(\overline{V}(x,y,z) - \varepsilon\right)^{+} - \left(V(x,z) - \varepsilon\right)^{+}\right)$$

$$+ \frac{\partial V}{\partial x}(x,z)g(x,\theta)s + s\left(\varphi'(x,y) - \frac{\partial k}{\partial x}(x,z)\Phi(x)\right)\theta \quad (4.36)$$

$$- s\frac{\partial k}{\partial x}(x,z)\left(f(x) + g(x,\theta)y\right) + s\left(\alpha'(x,y) - \frac{\partial k}{\partial x}(x,z)G(x)\right)d$$

$$+ a\frac{|d|^{2} + \left(\left(|\theta| - b - \lambda(e^{z})\right)^{+}\right)^{2}}{1 + \kappa(e^{z})} + s\left(h(x,y) + \tilde{g}(x,y,\theta)\overline{k}(x,y,z)\right)$$

Using (4.25), (4.31), (4.36) and the fact that $\left|v^{+} - w^{+}\right| \leq |v - w|$ for all $v, w \in \mathbb{R}$, we get:

$$\Lambda \leq -cV(x,z) + |s|\left|\frac{\partial k}{\partial z}(x,z)\right|\Gamma e^{-z}\left(\overline{V}(x,y,z) - \varepsilon\right)^{+} - M(x,y,z)s^{2}$$

$$+ \left|\frac{\partial V}{\partial z}(x,z)\right|\Gamma e^{-z}\left|\overline{V}(x,y,z) - V(x,z)\right| + a\frac{|d|^{2} + \left(\left(|\theta| - b - \lambda(e^{z})\right)^{+}\right)^{2}}{1 + \kappa(e^{z})}$$

$$+ |s|\left(|h(x,y)| + \left|\frac{\partial V}{\partial x}(x,z)\right||g(x,\theta)| + \left|\varphi'(x,y) - \frac{\partial k}{\partial x}(x,z)\Phi(x)\right||\theta|\right) \quad (4.37)$$

$$+ |s|\left(\left|\frac{\partial k}{\partial x}(x,z)\right||f(x) + g(x,\theta)y| + \left|\alpha'(x,y) - \frac{\partial k}{\partial x}(x,z)G(x)\right||d|\right)$$

Using the inequality

$$|s|\left|\alpha'(x,y) - \frac{\partial k}{\partial x}(x,z)G(x)\right||d|$$

$$\leq a\frac{|d|^{2}}{1 + \kappa(e^{z})} + \frac{1 + \kappa(e^{z})}{4a}\left|\alpha'(x,y) - \frac{\partial k}{\partial x}(x,z)G(x)\right|^{2}s^{2}$$

definition (4.27), the triangle inequality ($|f(x) + g(x,\theta)y| \leq |f(x)| + |g(x,\theta)||y|$) and inequality (4.26) we get from (4.37):

$$\Lambda \leq -cV(x,z) + |s|\left|\frac{\partial k}{\partial z}(x,z)\right|\Gamma e^{-z}\left(\bar{V}(x,y,z) - \varepsilon\right)^{+}$$

$$+ a\frac{2|d|^2 + \left(\left(|\theta| - b - \lambda(e^z)\right)^{+}\right)^2}{1 + \kappa(e^z)} + \frac{1}{2}s^2\left|\frac{\partial V}{\partial z}(x,z)\right|\Gamma e^{-z}$$

$$+ |s|\left(\left|\frac{\partial V}{\partial x}(x,z)\right|\mu(x) + \left|\frac{\partial k}{\partial x}(x,z)\right||f(x)| + \left|\frac{\partial k}{\partial x}(x,z)\right|\mu(x)|y| + |h(x,y)|\right) \quad (4.38)$$

$$+ |s|\left(\left|\varphi'(x,y) - \frac{\partial k}{\partial x}(x,z)\Phi(x)\right| + \left|\frac{\partial V}{\partial x}(x,z)\right|\mu(x) + \left|\frac{\partial k}{\partial x}(x,z)\right|\mu(x)|y|\right)|\theta|$$

$$- \left(M(x,y,z) - \frac{1 + \kappa(e^z)}{4a}\left|\alpha'(x,y) - \frac{\partial k}{\partial x}(x,z)G(x)\right|^2\right)s^2$$

The triangle inequalities $\left|\varphi'(x,y) - \frac{\partial k}{\partial x}(x,z)\Phi(x)\right| \leq \left|\frac{\partial k}{\partial x}(x,z)\Phi(x)\right| + |\varphi(x,y)|$ and $|y| \leq |s| + |k(x,z)|$ (recall that $s = y - k(x,z)$), inequalities (4.30), (4.38) and the fact that $\left(\bar{V}(x,y,z) - \varepsilon\right)^{+} \leq \bar{V}(x,y,z)$ give:

$$\Lambda \leq -cV(x,z) + |s|\left|\frac{\partial k}{\partial z}(x,z)\right|\Gamma e^{-z}\bar{V}(x,y,z) + a\frac{2|d|^2 + \left(\left(|\theta| - b - \lambda(e^z)\right)^{+}\right)^2}{1 + \kappa(e^z)}$$

$$+ \left(\left|\frac{\partial k}{\partial x}(x,z)\right|\mu(x) + \rho(x,y) + \frac{1}{2}\left|\frac{\partial V}{\partial z}(x,z)\right|\Gamma e^{-z}\right)s^2 + \left(\rho(x,y) + \left|\frac{\partial k}{\partial x}(x,z)\right|\mu(x)\right)s^2|\theta|$$

$$+ |s||x|\left(\left|\frac{\partial k}{\partial x}(x,z)\right|r(x) + \rho(x,y) + \left(\left|\frac{\partial k}{\partial x}(x,z)\right|\mu(x) + \mu(x) + \rho(x,y)\right)R(x,z)\right)$$

$$+ |s||x|\left(\rho(x,y) + R(x,z)\left(1 + \rho(x,y) + \mu(x) + \left|\frac{\partial k}{\partial x}(x,z)\right|\mu(x)\right)\right)|\theta|$$

$$- \left(M(x,y,z) - \frac{1 + \kappa(e^z)}{4a}\left|\alpha'(x,y) - \frac{\partial k}{\partial x}(x,z)G(x)\right|^2\right)s^2$$

$$(4.39)$$

Using the inequalities $\left|\frac{\partial k}{\partial x}(x,z)\right| \leq \frac{1}{2} + \frac{1}{2}\left|\frac{\partial k}{\partial x}(x,z)\right|^2$, $\left|\frac{\partial V}{\partial z}(x,z)\right| \leq \frac{1}{2} + \frac{1}{2}\left(\frac{\partial V}{\partial z}(x,z)\right)^2$, $\left|\frac{\partial k}{\partial z}(x,z)\right| \leq \frac{1}{2} + \frac{1}{2}\left(\frac{\partial k}{\partial z}(x,z)\right)^2$ and definition (4.27), we get from (4.39) and (4.34):

$$\Lambda \leq -cV(x,z) + \frac{\Gamma e^{-z}}{2}|s|\left(1+\left(\frac{\partial k}{\partial z}(x,z)\right)^2\right)V(x,z)$$

$$+\frac{\Gamma e^{-z}}{4}|s|^3\left(1+\left(\frac{\partial k}{\partial z}(x,z)\right)^2\right)+a\frac{2|d|^2+\left(\left(|\theta|-b-\lambda(e^z)\right)^+\right)^2}{1+\kappa(e^z)}$$

$$+\left(\frac{\mu(x)}{2}\left(1+\left|\frac{\partial k}{\partial x}(x,z)\right|^2\right)+\rho(x,y)+\frac{\Gamma e^{-z}}{4}\left(1+\left(\frac{\partial V}{\partial z}(x,z)\right)^2\right)\right)s^2 \quad (4.40)$$

$$+P(x,y,z)s^2|\theta|+|s||x|P(x,y,z)(1+|\theta|)$$

$$-\left(M(x,y,z)-\frac{1+\kappa(e^z)}{4a}\left|\alpha'(x,y)-\frac{\partial k}{\partial x}(x,z)G(x)\right|^2\right)s^2$$

The inequalities

$$|s| \leq \frac{1}{2}+\frac{1}{2}s^2$$

$$\frac{1}{2}|s|\left(1+\left(\frac{\partial k}{\partial z}(x,z)\right)^2\right)\Gamma e^{-z} \leq \frac{c}{4}+\frac{\Gamma^2 e^{-2z}}{4c}s^2\left(1+\left(\frac{\partial k}{\partial z}(x,z)\right)^2\right)^2$$

in conjunction with (4.40) give:

$$\Lambda \leq -\frac{3c}{4}V(x,z)+\frac{\Gamma^2 e^{-2z}}{4c}s^2\left(1+\left(\frac{\partial k}{\partial z}(x,z)\right)^2\right)^2 V(x,z)$$

$$+\frac{\Gamma e^{-z}}{8}(1+s^2)\left(1+\left(\frac{\partial k}{\partial z}(x,z)\right)^2\right)s^2 +a\frac{2|d|^2+\left(\left(|\theta|-b-\lambda(e^z)\right)^+\right)^2}{1+\kappa(e^z)}$$

$$+\left(\frac{\mu(x)}{2}\left(1+\left|\frac{\partial k}{\partial x}(x,z)\right|^2\right)+\rho(x,y)+\frac{\Gamma e^{-z}}{4}\left(1+\left(\frac{\partial V}{\partial z}(x,z)\right)^2\right)\right)s^2 \quad (4.41)$$

$$+P(x,y,z)s^2\left(|\theta|-b-\lambda(e^z)\right)+P(x,y,z)\left(b+\lambda(e^z)\right)s^2$$

$$+|s||x|P(x,y,z)\left(|\theta|-b-\lambda(e^z)\right)+|s||x|\left(1+b+\lambda(e^z)\right)P(x,y,z)$$

$$-\left(M(x,y,z)-\frac{1+\kappa(e^z)}{4a}\left|\alpha'(x,y)-\frac{\partial k}{\partial x}(x,z)G(x)\right|^2\right)s^2$$

Using the inequalities

$$|\theta|-b-\lambda(e^z) \leq \left(|\theta|-b-\lambda(e^z)\right)^+$$

$$|s||x|P(x,y,z)(b+1+\lambda(e^z))$$
$$\leq \frac{c}{4\sigma(x,z)}|x|^2 + \frac{\sigma(x,z)}{c}P^2(x,y,z)(b+1+\lambda(e^z))^2 s^2$$

we obtain from (4.41) and (4.33):

$$\Lambda \leq -\frac{3c}{4}V(x,z) + \frac{c}{4\sigma(x,z)}|x|^2 - \frac{c}{4}s^2$$
$$+a\frac{2|d|^2 + \left(\left(|\theta|-b-\lambda(e^z)\right)^+\right)^2}{1+\kappa(e^z)} + P(x,y,z)s^2\left(|\theta|-b-\lambda(e^z)\right)^+ \qquad (4.42)$$
$$+|s||x|P(x,y,z)\left(|\theta|-b-\lambda(e^z)\right)^+ - \frac{1+\kappa(e^z)}{2a}P^2(x,y,z)\left(s^2+|x|^2\right)s^2$$

Taking into account (4.23), (4.27), (4.42) and the inequalities

$$s^2 P(x,y,z)\left(|\theta|-b-\lambda(e^z)\right)^+$$
$$\leq \frac{a}{2}\frac{\left(\left(|\theta|-b-\lambda(e^z)\right)^+\right)^2}{1+\kappa(e^z)} + \frac{1+\kappa(e^z)}{2a}P^2(x,y,z)s^4$$

$$|s||x|P(x,y,z)\left(|\theta|-b-\lambda(e^z)\right)^+$$
$$\leq \frac{a}{2}\frac{\left(\left(|\theta|-b-\lambda(e^z)\right)^+\right)^2}{1+\kappa(e^z)} + \frac{1+\kappa(e^z)}{2a}|x|^2 P^2(x,y,z)s^2$$

we are in a position to conclude that (4.29) holds.

Finally, we proceed to show the validity of (4.28) with $\bar{\sigma}(x,y,z)$ defined by (4.32). Using the triangle inequality $|y| \leq |y-k(x,z)| + |k(x,z)|$ and the fact that $(|y-k(x,z)| + |k(x,z)|)^2 \leq 2|y-k(x,z)|^2 + 2|k(x,z)|^2$ we get:

$$|x|^2 + |y|^2 \leq |x|^2 + 2|y-k(x,z)|^2 + 2|k(x,z)|^2$$

The above inequality in conjunction with definition (4.27) and inequalities (4.30) gives:

$$|x|^2 + |y|^2 \leq \left(1+2R^2(x,z)\right)|x|^2 + 4\bar{V}(x,y,z) \qquad (4.43)$$

Using definition (4.27) and inequalities (4.23), (4.43) we obtain

$$|x|^2 + |y|^2 \leq \left(\left(1 + 2R^2(x,z)\right)\sigma(x,z) + 4\right)\bar{V}(x,y,z) \tag{4.44}$$

Combining (4.44) and (4.32) we obtain (4.28). The proof is complete. ◁

We next provide the proof of Theorem 1.

**Proof of Theorem 1:** By virtue of Theorem 2, it suffices to show that for every constants $b, \Gamma, \varepsilon, a, c > 0$ and for every $\kappa, \lambda \in K_\infty \cap C^\infty(\mathbb{R}_+)$ there exist functions $k, V \in C^\infty(\mathbb{R}^n \times \mathbb{R}^m \times \mathbb{R})$ with $V(0, 0, z) = k(0, 0, z) = 0$ for $z \in \mathbb{R}$, $V(x, y, z) > 0$ for $(x, y) \neq 0$, and a constant $M > 0$ such that:

(i) for each $R > 0$ the set $\{(x, y, z) \in \mathbb{R}^n \times \mathbb{R}^m \times [-R, R] : V(x, y, z) \leq R\}$ is bounded, and

(ii) the following inequalities hold:

$$|x|^2 + y_1^2 \leq MV(x, y, z), \text{ for all } (x, y, z) \in \mathbb{R}^n \times \mathbb{R}^m \times \mathbb{R} \tag{4.45}$$

$$\sum_{i=1}^{n-1} \frac{\partial V}{\partial x_i}(x,y,z)x_{i+1} + \frac{\partial V}{\partial x_n}(x,y,z)y_1 + \sum_{j=1}^{m-1} \frac{\partial V}{\partial y_j}(x,y,z)h_j(x,y_1,\ldots,y_j)$$

$$+ \sum_{j=1}^{m-1} \frac{\partial V}{\partial y_j}(x,y,z)\left(g_j(x,y_1,\ldots,y_j,\theta)y_{j+1} + \varphi'_j(x,y_1,\ldots,y_j)\theta + \alpha'_j(x,y_1,\ldots,y_j)d\right)$$

$$+ \frac{\partial V}{\partial y_m}(x,y,z)\left(h_m(x,y) + g_m(x,y,\theta)k(x,y,z) + \varphi'_m(x,y)\theta + \alpha'_m(x,y)d\right)$$

$$+ \frac{\partial V}{\partial z}(x,y,z)\Gamma e^{-z}\left(V(x,y,z) - \frac{\varepsilon^2}{2M}\right)^+ \leq -cV(x,y,z) + \frac{a}{M}\frac{|d|^2 + \left(\left(|\theta| - b - \lambda(e^z)\right)^+\right)^2}{1 + \kappa(e^z)}$$

for all $(x, y, z) \in \mathbb{R}^n \times \mathbb{R}^m \times \mathbb{R}$, $d \in \mathbb{R}^l$, $\theta \in \Theta$ \hfill (4.46)

Let arbitrary constants $b, \Gamma, \varepsilon, a, c > 0$ and arbitrary functions $\kappa, \lambda \in K_\infty \cap C^\infty(\mathbb{R}_+)$ be given. The functions $k, V \in C^\infty(\mathbb{R}^n \times \mathbb{R}^m \times \mathbb{R})$ are constructed inductively by using Lemma 4.

Define the matrices

$$A = \{a_{i,j} : i, j = 1, \ldots, n\}, \; b = \begin{bmatrix} 0 \\ \vdots \\ 0 \\ 1 \end{bmatrix} \tag{4.47}$$

where

$$a_{i,j} = 1 \text{ for } j = i+1, \; i = 1, \ldots, n-1 \text{ and } a_{i,j} = 0 \text{ if otherwise} \tag{4.48}$$

Notice that equations (2.1) and definitions (4.47), (4.48) and the fact that $x_{n+1} = y_1$ imply the equation:

$$\dot{x} = Ax + by_1 \tag{4.49}$$

Since the pair of matrices $(A,b)$ is a controllable pair there exists a vector $\omega \in \mathbb{R}^n$ and a positive definite, symmetric matrix $P \in \mathbb{R}^{n \times n}$ such that

$$P(A + b\omega') + (A + b\omega')'P \leq -2^m cP \tag{4.50}$$

Next, we notice that the functions

$$V_1(x, y_1, z) := x'Px + \frac{1}{2}(y_1 - \omega'x)^2, \quad k_1(x, y_1, z) := -\frac{G(x, y_1, z)}{\eta_1(x, y_1)}(y_1 - \omega'x) \tag{4.51}$$

with

$$G(x, y_1, z) := \frac{M}{2^{m+1}c}\left(K + r(x, y_1)(1 + |\omega|)(1 + b + \lambda(e^z))\right)^2 + |\omega'b| + 2^{m-2}c$$
$$+ \frac{M(1 + \kappa(e^z))}{2^{3-m}a}\left(|\alpha_1(x, y_1)|^2 + r^2(x, y_1)(|x|^2 + y_1^2)\right) + r(x, y_1)(1 + b + \lambda(e^z)) \tag{4.52}$$

where $K = |2b'P - \omega'A - \omega'b\omega'|$, $\eta_1$ is the function involved in (2.3) and $r : \mathbb{R}^n \times \mathbb{R} \to (0, +\infty)$ is a smooth function that satisfies

$$|h_1(x, y_1)| + |\varphi_1(x, y_1)| \leq r(x, y_1)(|x| + |y_1|), \text{ for all } (x, y_1) \in \mathbb{R}^n \times \mathbb{R} \tag{4.53}$$

satisfy the following inequality for all $(x, y_1, z) \in \mathbb{R}^n \times \mathbb{R} \times \mathbb{R}$, $d \in \mathbb{R}^l$, $\theta \in \Theta$

$$\sum_{i=1}^{n-1} \frac{\partial V_1}{\partial x_i}(x, y_1, z)x_{i+1} + \frac{\partial V_1}{\partial x_n}(x, y_1, z)y_1$$
$$+ \frac{\partial V_1}{\partial y_1}(x, y_1, z)\left(h_1(x, y_1) + g_1(x, y_1, \theta)k_1(x, y_1, z) + \varphi_1'(x, y_1)\theta + \alpha_1'(x, y_1)d\right)$$
$$+ \frac{\partial V_1}{\partial z}(x, y_1, z)\Gamma e^{-z}\left(V(x, y_1, z) - \frac{\varepsilon^2}{2M}\right)^+ \tag{4.54}$$
$$\leq -2^{m-1}cV_1(x, y_1, z) + 2^{1-m}\frac{a}{M}\frac{|d|^2 + \left((|\theta| - b - \lambda(e^z))^+\right)^2}{1 + \kappa(e^z)}$$

where $M > 0$ is any constant for which the following inequality holds:

$$M\left(x'Px + \frac{1}{2}(y_1 - \omega'x)^2\right) \geq |x|^2 + y_1^2, \text{ for all } (x, y_1) \in \mathbb{R}^n \times \mathbb{R} \tag{4.55}$$

The existence of a constant $M > 0$ for which (4.55) holds follows from definition (4.51) and the fact that $x'Px + \frac{1}{2}(y_1 - \omega'x)^2$ is a positive definite quadratic form on $\mathbb{R}^{n+1}$.

Inequality (4.54) is a direct consequence of definitions (4.51), (4.52), inequalities (2.3), (4.50), (4.53) and (4.55) (which implies the inequality $M\,x'Px \geq |x|^2$ for all $x \in \mathbb{R}^n$) as well as the inequalities

$$|x|+|y_1| \leq (1+|\omega|)|x|+|y_1-\omega'x|$$

$$|y_1-\omega'x||\alpha_1(x,y_1)||d|$$
$$\leq \frac{2^{1-m}a|d|^2}{M(1+\kappa(e^z))} + \frac{M(1+\kappa(e^z))}{2^{3-m}a}|\alpha_1(x,y_1)|^2(y_1-\omega'x)^2$$

$$|y_1-\omega'x|r(x,y_1)(|x|+|y_1|)(|\theta|-b-\lambda(e^z))^+ \leq 2^{1-m}a\frac{\left((|\theta|-b-\lambda(e^z))^+\right)^2}{M(1+\kappa(e^z))}$$
$$+\frac{M(1+\kappa(e^z))}{2^{3-m}a}r^2(x,y_1)(|x|+|y_1|)^2(y_1-\omega'x)^2$$

$$(|x|+|y_1|)^2 \leq 2|x|^2 + 2y_1^2$$

$$|\theta|-b-\lambda(e^z) \leq (|\theta|-b-\lambda(e^z))^+$$

$$|y_1-\omega'x||x|\left(K+r(x,y_1)(1+|\omega|)(1+b+\lambda(e^z))\right) \leq 2^{m-1}c\frac{|x|^2}{M}$$
$$+\frac{M}{2^{m+1}c}\left(K+r(x,y_1)(1+|\omega|)(1+b+\lambda(e^z))\right)^2(y_1-\omega'x)^2$$

which give inequality (4.54).

Applying Lemma 4 exactly $m-1$ times, we construct functions $k_j \in C^\infty(\mathbb{R}^n \times \mathbb{R}^j \times \mathbb{R})$ for $j=2,...,m$ with $k_j(0,0,z)=0$ for $z \in \mathbb{R}$, such that (4.46) holds for the function

$$V(x,y,z) = x'Px + \frac{1}{2}(y_1-\omega'x)^2 + \frac{1}{2}\sum_{j=2}^{m}\left(y_j - k_{j-1}(x,y_1,...,y_{j-1},z)\right)^2 \qquad (4.56)$$

Notice that inequality (4.55) and (4.56) guarantee the validity of (4.45).

The proof is complete.    ◁

Finally, we provide the proof of Theorem 3.

**Proof of Theorem 3:** By virtue of Theorem 2, it suffices to show that for every constants $b, \Gamma, \varepsilon, a, c > 0$ and for every $\kappa, \lambda \in K_\infty \cap C^\infty(\mathbb{R}_+)$ there exist functions $k, V \in C^\infty(\mathbb{R}^n \times \mathbb{R})$ with $V(0, z) = k(0, z) = 0$ for $z \in \mathbb{R}$, $V(x, z) > 0$ for $x \neq 0$, such that:

(i) for each $R > 0$ the set $\{(x, z) \in \mathbb{R}^n \times [-R, R] : V(x, z) \leq R\}$ is bounded, and

(ii) the following inequalities hold:

$$\frac{1}{2} x_1^2 \leq V(x, z), \text{ for all } (x, z) \in \mathbb{R}^n \times \mathbb{R} \qquad (4.57)$$

$$\sum_{i=1}^{n-1} \frac{\partial V}{\partial x_i}(x, z)\left(h_i(x_1, ..., x_i) + g_i(x_1, ..., x_i, \theta)x_{i+1} + \varphi_i'(x_1, ..., x_i)\theta + \alpha_i'(x_1, ..., x_i)d\right)$$

$$+ \frac{\partial V}{\partial x_n}(x, z)\left(h_n(x) + g_n(x, \theta)k(x, z) + \varphi_n'(x)\theta + \alpha_n'(x)d\right)$$

$$+ \frac{\partial V}{\partial z}(x, z)\Gamma e^{-z}\left(V(x, z) - \frac{\varepsilon^2}{2}\right)^+ \leq -cV(x, z) + a\frac{|d|^2 + \left(\left(|\theta| - b - \lambda(e^z)\right)^+\right)^2}{1 + \kappa(e^z)}$$

$$\text{for all } (x, z) \in \mathbb{R}^n \times \mathbb{R}, \, d \in \mathbb{R}^l, \, \theta \in \Theta \qquad (4.58)$$

Let arbitrary constants $b, \Gamma, \varepsilon, a, c > 0$ and arbitrary functions $\kappa, \lambda \in K_\infty \cap C^\infty(\mathbb{R}_+)$ be given. The functions $k, V \in C^\infty(\mathbb{R}^n \times \mathbb{R})$ are constructed inductively by using Lemma 4.

Indeed, we notice that the functions

$$V_1(x_1, z) := \frac{1}{2} x_1^2, \quad k_1(x_1, z) := -\frac{M(x_1, z)}{\eta_1(x_1)} x_1 \qquad (4.59)$$

with

$$M(x_1, z) := \left(b + 1 + \lambda(e^z)\right)r(x_1) + \frac{1 + \kappa(e^z)}{2^{3-n}a}\left(|\alpha_1(x_1)|^2 + r^2(x_1)x_1^2\right) + 2^{n-2}c \qquad (4.60)$$

where $\eta_1$ is the function involved in (2.19) and $r : \mathbb{R} \to (0, +\infty)$ is a smooth function that satisfies

$$|h_1(x_1)| + |\varphi_1(x_1)| \leq r(x_1)|x_1|, \text{ for all } x_1 \in \mathbb{R} \qquad (4.61)$$

satisfy the following inequality for all $(x_1, z) \in \mathbb{R} \times \mathbb{R}$, $d \in \mathbb{R}^l$, $\theta \in \Theta$

$$\frac{\partial V_1}{\partial x_1}(x_1, z)\left(h_1(x_1) + g_1(x_1, \theta)k_1(x_1, z) + \varphi_1'(x_1)\theta + \alpha_1'(x_1)d\right)$$

$$+ \frac{\partial V_1}{\partial z}(x_1, z)\Gamma e^{-z}\left(V(x_1, z) - \frac{\varepsilon^2}{2}\right)^+ \qquad (4.62)$$

$$\leq -2^{n-1}cV_1(x_1, z) + 2^{1-n}a\frac{|d|^2 + \left(\left(|\theta| - b - \lambda(e^z)\right)^+\right)^2}{1 + \kappa(e^z)}$$

Inequality (4.62) is a direct consequence of definitions (4.59), (4.60), inequalities (2.19) and (4.61) as well as the inequalities

$$|x_1||\alpha_1(x_1)||d| \leq 2^{1-n} a \frac{|d|^2}{1+\kappa(e^z)} + \frac{1+\kappa(e^z)}{2^{3-n} a} |\alpha_1(x_1)|^2 x_1^2$$

$$|\theta| - b - \lambda(e^z) \leq (|\theta| - b - \lambda(e^z))^+$$

$$r(x_1)x_1^2 (|\theta| - b - \lambda(e^z))^+$$

$$\leq \frac{1+\kappa(e^z)}{2^{3-n} a} r^2(x_1) x_1^4 + 2^{1-n} a \frac{((|\theta| - b - \lambda(e^z))^+)^2}{1+\kappa(e^z)}$$

which give inequality (4.62).

Applying Lemma 4 exactly $n-1$ times, we construct functions $k_i \in C^\infty(\mathbb{R}^i \times \mathbb{R})$ for $i = 2,...,n$ with $k_i(0,z) = 0$ for $z \in \mathbb{R}$, such that (4.58) holds for the function

$$V(x,z) = \frac{1}{2} x_1^2 + \frac{1}{2} \sum_{i=2}^n (x_i - k_{i-1}(x_1,...,x_{i-1},z))^2 \qquad (4.63)$$

Notice that equality (4.63) guarantees the validity of (4.57).

The proof is complete. ◁

## 5. Concluding Remarks

The DADS controller is a direct adaptive control scheme that can be applied to systems that satisfy a matching condition (see [12]) or systems in parametric strict feedback form. It is clear that further research is needed in order to identify more classes of nonlinear systems to which the DADS controller can be applied. On the other hand, the examples in [12] show that the zero p-OAG property $\limsup_{t \to +\infty}(|Y(t)|) \leq \alpha$ with a constant $\alpha > 0$ independent of $\theta$ and for an appropriate output map $Y = h(x)$ cannot be achieved by feedback control (static or dynamic) for an arbitrary nonlinear control system: there are controllable nonlinear systems which defy this property. Therefore, the identification of more classes of nonlinear systems to which the DADS controller can be applied is not an easy problem.

For indirect adaptive schemes, apart from regulation, identification of the constant unknown parameters $\theta$ is also achieved in certain cases. The development of novel (and robust with respect to disturbances like $d$) identifiers for direct adaptive schemes (like DADS) is an interesting problem that requires special attention.